\documentclass[11pt]{article}

\title{L\'evy-driven Volterra equations in space and time}

\author{
Carsten Chong\thanks{Center for Mathematical Sciences, Technische Universit\"at M\"unchen, Boltzmannstra\ss e 3, 85748 Garching, Germany, e-mail: carsten.chong@tum.de, url: www.statistics.ma.tum.de}
}

\usepackage{latexsym}
\usepackage{amssymb}
\usepackage{amsmath}
\usepackage{amsthm}
\usepackage{mathrsfs}
\usepackage[numbers]{natbib}
\usepackage{stmaryrd}
\usepackage{url}
\usepackage{longtable}
\usepackage{color}
\usepackage{dsfont}
\usepackage[T1]{fontenc}
\usepackage{lmodern}

\usepackage[latin1]{inputenc}
\usepackage{enumitem}
\setenumerate{leftmargin=*}
\usepackage{amsopn}
\usepackage{comment}

\newcommand{\bfi}{\begin{fig}}
\newcommand{\efi}{\end{fig}}
\newcommand{\btab}{\begin{tab}}
\newcommand{\etab}{\end{tab}}
\newcommand{\barr}{\begin{array}}
\newcommand{\earr}{\end{array}}
\newcommand{\beq}{\begin{equation}}
\newcommand{\eeq}{\end{equation}}
\newcommand{\bdis}{\begin{displaymath}}
\newcommand{\edis}{\end{displaymath}\noindent}

\newcommand{\bbn}{\mathbb{N}}

\newcommand{\bbr}{\mathbb{R}}

\newcommand{\bbe}{\mathbb{E}}

\newcommand{\bbp}{\mathbb{P}}
\newcommand{\bbb}{\mathbb{B}}
\newcommand{\bbf}{\mathbb{F}}

\newcommand{\bone}{\mathds 1}

\newcommand{\eqd}{\stackrel{d}{=}}

\newcommand{\halmos}{\quad\hfill $\Box$}

\newcommand{\cals}{{\cal S}}

\newcommand{\calf}{{\cal F}}

\newcommand{\cale}{{\cal E}}
\newcommand{\calp}{{\cal P}}
\newcommand{\calb}{{\cal B}}

\newcommand{\calm}{{\cal M}}

\newcommand{\al}{{\alpha}}
\newcommand{\la}{{\lambda}}
\newcommand{\La}{{\Lambda}}
\newcommand{\eps}{{\epsilon}}

\newcommand{\ga}{{\gamma}}

\newcommand{\si}{{\sigma}}
\newcommand{\om}{{\omega}}
\newcommand{\Om}{{\Omega}}

\newcommand{\ov}{\overline}
\newcommand{\un}{\underline}

\newcommand{\dd}{\mathrm{d}}
\newcommand{\ee}{\mathrm{e}}

\newcommand{\bb}{\mathrm{b}}
\newcommand{\cc}{\mathrm{c}}

\newcommand{\Comp}{\mathrm{c}}
\newcommand{\loc}{\mathrm{loc}}
\newcommand{\diam}{\mathrm{diam}}

\newcommand{\ppp}{^\mathrm{p}\!}

\newcommand{\pf}{\mathfrak{p}}
\newcommand{\qf}{\mathfrak{q}}

\makeatletter
\newcommand{\opnorm}{\@ifstar\@opnorms\@opnorm}
\newcommand{\@opnorms}[1]{%
  \left|\mkern-1.5mu\left|\mkern-1.5mu\left|
   #1
  \right|\mkern-1.5mu\right|\mkern-1.5mu\right|
}
\newcommand{\@opnorm}[2][]{%
  \mathopen{#1|\mkern-1.5mu#1|\mkern-1.5mu#1|}
  #2
  \mathclose{#1|\mkern-1.5mu#1|\mkern-1.5mu#1|}
}
\makeatother

\newtheoremstyle{neu}
    {11pt}      
    {11pt}      
    {}                  
    {}          
    {\bfseries} 
    {}          
    {1em}  
    {\textbf{\thmname{#1}\thmnumber{ #2}\thmnote{ (#3)}}}          
    
\newtheoremstyle{proof}
    {11pt}      
    {11pt}      
    {}                  
    {}          
    {\bfseries} 
    {}            
    {1em}          
    {\textbf{\thmname{#1}.}}          

\newtheorem{Theorem}{Theorem}[section]
\newtheorem{Corollary}[Theorem]{Corollary}
\newtheorem{Lemma}[Theorem]{Lemma}
\newtheorem{Proposition}[Theorem]{Proposition}
\theoremstyle{neu}
\newtheorem{Definition}[Theorem]{Definition}
\newtheorem{Example}[Theorem]{Example}
\newtheorem{Remark}[Theorem]{Remark}
\newtheorem{Assumption}{Assumption}

\theoremstyle{proof}
\newtheorem{Proof}{Proof}

\newcommand{\bthm}{\begin{Theorem}}
\newcommand{\ethm}{\end{Theorem}}

\newcommand{\bcor}{\begin{Corollary}}
\newcommand{\ecor}{\end{Corollary}}

\newcommand{\blem}{\begin{Lemma}}
\newcommand{\elem}{\end{Lemma}}

\newcommand{\bprop}{\begin{Proposition}}
\newcommand{\eprop}{\end{Proposition}}

\newcommand{\bdf}{\begin{Definition}}
\newcommand{\edf}{\halmos \end{Definition}}

\newcommand{\bex}{\begin{Example}}
\newcommand{\eex}{\halmos \end{Example}}

\newcommand{\brem}{\begin{Remark}}
\newcommand{\erem}{\halmos \end{Remark}}

\newcommand{\bass}{\begin{Assumption}}
\newcommand{\eass}{\halmos \end{Assumption}}

\newcommand{\bpr}{\begin{Proof}}
\newcommand{\epr}{\halmos \end{Proof}}

\newcommand{\benu}{\begin{enumerate}}
\newcommand{\eenu}{\end{enumerate}}

\newcommand{\bit}{\begin{itemize}}
\newcommand{\eit}{\end{itemize}}

\newcommand{\bff}{\textbf}

\numberwithin{equation}{section}

\begin{document}


\maketitle

\begin{abstract}
We investigate nonlinear stochastic Volterra equations in space and time that are driven by L\'evy bases. Under a Lipschitz condition on the nonlinear term, we give existence and uniqueness criteria in weighted function spaces that depend on integrability properties of the kernel and the characteristics of the L\'evy basis. Particular attention is devoted to equations with stationary solutions, or more generally, to equations with infinite memory, that is, where the time domain of integration starts at minus infinity. Here, in contrast to the case where time is positive, the usual integrability conditions on the kernel are no longer sufficient for the existence and uniqueness of solutions, but we have to impose additional size conditions on the kernel and the L\'evy characteristics. Furthermore, once the existence of a solution is guaranteed, we analyse its asymptotic stability, that is, whether its moments remain bounded when time goes to infinity. Stability is proved whenever kernel and characteristics are small enough, or the nonlinearity of the equation exhibits a fractional growth of order strictly smaller than one. The results are applied to the stochastic heat equation for illustration.
\end{abstract}

\vfill

\noindent
\begin{tabbing}
\> secondary: \,\,\,60G51, 60G57
\end{tabbing}

\vspace{1cm}


\vspace{0.5cm}

\newpage

\section{Introduction}
In this paper we investigate stochastic tempo--spatial Volterra equations of the following form:
\beq\label{SPDE-var} Y(t,x)=Y_0(t,x)+\int_I \int_{\bbr^d} G(t,x;s,y)\si(Y(s,y))\,\La(\dd s,\dd y),\quad (t,x)\in I\times\bbr^d. \eeq
Here, $Y_0$ is a given stochastic process, $I$ is a real time interval, $G$ a deterministic kernel function and $\si$ a deterministic function. Apart from $Y_0$, the stochasticity of \eqref{SPDE-var} comes from its integrator $\La$, which is an infinitely divisible independently scattered random measure, or a \emph{L\'evy basis} for short. 

While the theory of deterministic Volterra equations is very well studied by now (see, for example, the monograph \citep{Gripenberg90}), the literature on Volterra equations with stochastic integrators is considerably smaller. If no space is involved, \citep{Protter85} proves existence and uniqueness for general semimartingale integrators 
under differentiability assumptions on the kernel $G$. In the special case of L\'evy-driven stochastic delay equations, the asymptotic behaviour of solutions and the existence of stationary solutions are discussed in \citep{Reiss06}. As soon as the kernel becomes explosive, existence and uniqueness results have been found for Brownian integrators, see \citep{Couchran95, Coutin01, Wang08}. In the tempo--spatial case, singular kernels are typically encountered in the theory of stochastic PDEs, with two main approaches having become established in this context: on the one hand, there is the functional analytic approach that treats infinite-dimensional stochastic evolution equations as ordinary SDEs with irregular coefficients driven by Hilbert or Banach space-valued L\'evy processes; see, for instance, \citep{Peszat07} for an excellent account on this subject; or see the recent paper \citep{Kovacs15} for the treatment of Volterra-type equations within this framework. On the other hand, there is the random field approach that directly considers \eqref{SPDE-var} as a scalar-valued equation driven by a multi-parameter L\'evy noise. In the Gaussian case, the two approaches have been compared in \citep{Dalang11}, in the general L\'evy case, this problem seems to be open.

Since our treatment of \eqref{SPDE-var} will be within the random field approach, we review the existing literature in this field in more detail: based on the seminal work \citep{Walsh86}, which uses equations of type \eqref{SPDE-var} in order to solve certain stochastic PDEs driven by Gaussian white noise, several attempts have been made to generalize Walsh's method to other noise types. One possibility is, for instance, to consider Gaussian noise that is white in time but coloured in space, which is proposed in \citep{Dalang99}. Leaving the Gaussian world, \citep{Albeverio98,Applebaum00} study the stochastic heat equation driven by L\'evy white noise. However, since both references still employ the $L^2$-theory of Walsh, they are confronted with the uncomfortable fact that the stochastic heat equation will have no solutions in dimensions greater than $1$, cf. \citep[pp.~328ff.]{Walsh86}. This is due to the bad integrability properties of the heat kernel that plays the role of $G$ in \eqref{SPDE-var}: it 
is square-integrable only for $d=1$.

Therefore, the passage from the $L^2$- to an $L^p$-framework, $p\in(0,2]$, is inevitable. The first paper that discusses L\'evy-driven stochastic PDEs in an $L^p$-framework with $p\in[1,2]$ is, to our best knowledge, \citep{SLB98}. Under the usual Lipschitz condition on $\si$, existence and uniqueness for \eqref{SPDE-var} are proved when $G$ is the heat kernel and $\La$ a homogeneous L\'evy basis that is either a martingale measure or of locally finite variation. In \citep{Mueller98, Mytnik02} a specific equation that goes beyond the results of \citep{SLB98} is studied: they take the non-Lipschitz coefficient $\si(x)=x^\beta$ with $\beta\neq1$ and an $\al$-stable spectrally positive L\'evy basis for $\La$, where $\al\in(0,1)$ and $\al\in(1,2)$, respectively. Finally, \citep{Balan14} treats the Lipschitz case with $\al$-stable $\La$ where $\al\neq1$. In all articles mentioned so far, the time horizon is $I=\bbr_+$.

Let us also point out that processes of the form \eqref{SPDE-var} are closely related to a class of random fields that are called \emph{ambit processes} and have found applications in physics, finance, biology among other disciplines; see \citep{BN04, BN11-2, BN15, Podolskij14} for more details. This class of processes takes the form
\beq\label{ambit} Y(t,x)= \mu + \int_{A(t,x)} G(t,x;s,y)\si(s,y)\,\La(\dd s,\dd y),\quad (t,x)\in \bbr\times\bbr^d,\eeq
where $A(t,x)$, the so-called ambit sets, are certain subsets of $\bbr\times\bbr^d$, $\mu\in\bbr$ is a constant and $\si$ is some given random field. As we can see, the major difference to \eqref{SPDE-var} is that the random field $\si$ in \eqref{ambit} is given by a function of $Y$ in \eqref{SPDE-var}. Once a solution to \eqref{SPDE-var} is found, it is a special type of ambit processes. For the connection between ambit processes and stochastic PDEs, we refer to \citep{BN11-2}.

The paper is organized as follows: after we have provided all necessary background information in Section~\ref{prelim}, we start to discuss \eqref{SPDE-var} in Section~\ref{localex} for $I=\bbr_+$. In Theorem~\ref{main1} we establish existence and uniqueness conditions for \eqref{SPDE-var} in $L^p$-spaces for $p\in(0,2]$ under Lipschitz conditions on $\si$. They generalize the results mentioned in the literature review to kernels $G$ that need not be of convolution type or related to stochastic PDEs, as well as to L\'evy bases that are combinations of martingale and finite variation parts, and whose characteristics are potentially inhomogeneous in space and time. The most stringent condition in Theorem~\ref{main1} is that, loosely speaking, $\La$ must have a moment structure that is at least as nice as its variation structure. This, for instance, a priori excludes any stable L\'evy basis. An extension to such cases is provided in Theorem~\ref{main2} if $\La$ only has finitely many large jumps on finite time intervals. Using localization methods as in \citep{Balan14}, we are able to reduce the situation to the framework of Theorem~\ref{main1} and prove existence and uniqueness of solutions this way. Beyond that, if $\si$ has sublinear growth, we prove that they have finite $L^p$-moments for some $p\in(0,2]$. 

In Section~\ref{globalex}, we extend the results from Section~\ref{localex} to the case of infinite memory, which, to our knowledge, has not been considered before in the literature. More precisely, we investigate existence and uniqueness for \eqref{SPDE-var} when $I=\bbr$ (Theorem~\ref{main3}), which turns out to be much more involved than the case $I=[0,\infty)$. First, the method of Theorem~\ref{main2} will no longer work, that is, $\La$ is required to have a good moment structure. Second, and more importantly, an explicit size condition on $G$, $\si$ and $\La$ comes into play, which is already a characteristic feature of deterministic Volterra equations, see Example~\ref{exexpfunc}. Therefore, detailed $L^p$-estimates for the stochastic integral in \eqref{SPDE-var} are required. Furthermore, under certain conditions on $Y_0$, one can improve the results by using weighted $L^p$-spaces. If $G$ is a kernel of convolution form and $\La$ is homogeneous in space and time, the stationarity of the solution is discussed in Theorem~\ref{statsol}. Section~\ref{globalex} is round off with some results concerning the $L^p$-continuity of the solution $Y$ and its continuous dependence on $Y_0$; see Theorem~\ref{further}.

In Section~\ref{assst} we assume that we have already found a solution to \eqref{SPDE-var} that is $L^p$-bounded up to time $T$ for every $T\in \bbr_+$. We want to address the question when the solution remains $L^p$-bounded as $T\to\infty$. An affirmative answer is given under two types of conditions (Theorem~\ref{assstat}): first, if $G$, $\si$ and $\La$ are small enough, a feature that we have already encountered in Theorem~\ref{main3} and that is also similar to the conditions in \citep{Reiss06} in the context of stationary solutions to stochastic delay equations; and second, if the function $\si$ is of sublinear growth. Both conditions are intrinsic for Volterra-type equations as a deterministic example shows, see Example~\ref{inftyex}. 

In Sections~\ref{localex} to \ref{assst}, we illustrate all our results by means of the stochastic heat equation, see Examples~\ref{ex1}, \ref{ex2}, \ref{ex3} and \ref{ex4}.

Finally, Section~\ref{proofs} contains several lemmata needed for the proof of the main theorems, which is carried out in Section~\ref{proofmain}.

\section{Preliminaries}\label{prelim}

We begin with a table of frequently used notations and abbreviations:
\begin{longtable}{p{2 cm} p{13 cm}}
$\bbr_+$ & the set $[0,\infty)$ of \emph{positive} real numbers, while \emph{strict positivity} excludes $0$;\\
$\bar\bbr$ & the extended real line $\bbr\cup\{\pm \infty\}$;\\
$\bbn$ & the set $\{1, 2, \ldots\}$ of natural numbers;\\
$I$ & either $I=\bbr_+$ or $I=\bbr$;\\
$I_T$ & $I\cap(-\infty,T]$ for some $T\in\bbr\cup\{\infty\}$;\\
$p^\ast$ & $p\vee 1$ for $p\in[0,\infty)$;\\
$|z|^r_s$ & $|z|^r\bone_{\{|z|>1\}} + |z|^s \bone_{\{|z|\leq1\}}$ for $r,s,z \in\bbr$;\\
$\bbb$ & a stochastic basis $(\Om,\calf,\bbf=(\calf_t)_{t\in I},\bbp)$ satisfying the usual hypotheses of right-continuity and completeness that is large enough to support all random elements of this paper;\\
$\tilde\Om$ & $\tilde\Om := \Om\times I \times\bbr^d$ for some $d\in\bbn\cup\{0\}$ with the convention $\bbr^0:=\{1\}$;\\
$\tilde\calp$ & depending on the context, either the \emph{tempo--spatial predictable $\si$-field} $\calp\otimes\calb(\bbr^d)$ where $\calp$ is the usual predictable $\si$-field and $\calb(\bbr^d)$ is the Borel $\si$-field on $\bbr^d$, or the class of \emph{predictable} (i.e. $\tilde\calp$-measurable) mappings $\tilde\Om\to\bar\bbr$;\\
$\tilde\calp_\bb$ & the collection of all sets $A\in\tilde\calp$ such that there exists $k\in\bbn$ with $A\subseteq \Om\times(I\cap[-k,k])\times[-k,k]^d$;\\
$\calb_\bb$ & the collection of all bounded Borel sets in $I\times\bbr^d$;\\
$\llbracket R,S \rrbracket$ & $\{(\om,t)\in\Om\times I\colon R(\om)\leq t\leq S(\om)\}$ for two $\bbf$-stopping times $R,S$, analogously for the other stochastic intervals;\\
$|\mu|$ & the total variation measure of a signed Borel measure $\mu$;\\
$x+A$ & $\{x+a\colon a\in A\}$ for $x\in\bbr^d$ and $A\subseteq \bbr^d$;\\
$A^\Comp$ & $\bbr^d\setminus A$ for $A\subseteq \bbr^d$;\\
$\diam(A)$ & $\sup \{|x-y|\colon x,y\in A\}$ for $A\subseteq \bbr^d$;\\
$(x,y]$ & $\{z\in\bbr^d\colon x_i<z_i\leq y_i \text{ for all } i=1,\ldots,d\}$ for $x,y\in\bbr^d$;\\
$L^p$ & the usual spaces $L^p(\Om,\calf,\bbp)$ for $p\in[0,\infty)$ endowed with the topologies induced by $\|X\|_{L^p} := \bbe[|X|^p]^{1/p^\ast}$ for $p\in(0,\infty)$ and $\|X\|_{L^0} := \bbe[|X|\wedge 1]$ for $p=0$;
\end{longtable}

In model \eqref{SPDE-var}, $\La$ will always be a \emph{L\'evy basis} on $I\times\bbr^d$, that is, a mapping $\La\colon \tilde\calp_\bb \to L^0$ with the following properties:
\benu
  \item $\La(\emptyset)=0$ a.s.
\item For every sequence $(A_i)_{i\in\bbn}$ of pairwise disjoint sets in $\tilde\calp_\bb$ with $\bigcup_{i=1}^\infty A_i\in\tilde\calp_\bb$ we have
\[ \La\left(\bigcup_{i=1}^\infty A_i\right)=\sum_{i=1}^\infty \La(A_i)\quad\text{in } L^0.\]
\item For all $A\in\tilde\calp_\bb$ with $A\subseteq \Om\times I_t \times\bbr^d$ for some $t\in I$, the random variable $\La(A)$ is $\calf_t$-measurable.
\item For all $A\in\tilde\calp_\bb$, $t\in I$ and $\Om_0\in\calf_t$, we have
\[\La\big(A\cap(\Om_0\times (t,\infty)\times \bbr^d)\big)=\bone_{\Om_0} \La\big(A\cap(\Om\times(t,\infty)\times \bbr^d)\big)\quad\text{a.s.}\]
\item If $(B_i)_{i\in\bbn}$ is a sequence of pairwise disjoint sets in $\calb_\bb$, then $(\La(\Om\times B_i))_{i\in\bbn}$ is a sequence of independent random variables. Furthermore, if $B\in\calb_\bb$ satisfies $B\subseteq (t,\infty)\times\bbr^d$ for some $t\in I$, then $\La(\Om\times B)$ is independent of $\calf_t$.
\item For all $B\in\calb_\bb$, $\La(\Om\times B)$ has an infinitely divisible distribution.
\item For all $t\in I$ and $k\in\bbn$ we have $\La(\Om\times\{t\}\times[-k,k]^d)=0$ a.s. 
\eenu

L\'evy bases are originally called infinitely divisible independently scattered random measures in \citep{Rajput89}; the short terminology has been introduced in \citep{BN04}. Conditions (3) and (4) are added to ensure that L\'evy bases are ``adapted'' to the underlying stochastic bases, see e.g. \citep{Chong14}. Just as L\'evy processes are semimartingales in the purely temporal case, L\'evy bases are random measures, that is, stochastic integrators in space--time. In other words, it is possible to develop an It\^o stochastic integration theory for L\'evy bases. Let us briefly recall this; all details can be found in \citep[Chap.~3]{Bichteler02} and \citep{Bichteler83}. Starting with simple integrands $H\in\cals$, that is, $H=\sum_{i=1}^r a_i\bone_{A_i}$ with $r\in\bbn$, real numbers $a_i$ and sets $A_i\in\tilde\calp_\bb$, we define the stochastic integral in the canonical way:
\[ \int_I\int_{\bbr^d} H(t,x)\,\La(\dd t,\dd x) := \sum_{i=1}^r a_i \La(A_i). \]
Given a general predictable function $H\in\tilde\calp$, we introduce the \emph{Daniell mean}
\[ \|H\|_\La := \sup_{S\in\cals, |S|\leq|H|} \left\|\int_I\int_{\bbr^d} S(t,x)\,\La(\dd t,\dd x)\right\|_{L^0}, \]
and define the class of \emph{integrable} functions $L^0(\La)$ as the closure of $\cals$ under the Daniell mean $\|\cdot\|_\La$. This is to say that $H\in\tilde\calp$ is integrable w.r.t. $\La$ if and only if there exists a sequence $(S_n)_{n\in\bbn}$ of elements in $\cals$ such that $\|H-S_n\|_\La\to0$ as $n\to\infty$. Then the \emph{stochastic integral}
\[ \int_I \int_{\bbr^d} H(t,x)\,\La(\dd t,\dd x) := \lim_{n\to\infty} \int_I\int_{\bbr^d} S_n(t,x)\,\La(\dd t,\dd x) \]
as a limit in probability exists and does not depend on the chosen sequence $(S_n)_{n\in\bbn}$. Moreover, defining 
\[ H\cdot\La_t := \int_{I_t}\int_{\bbr^d} H(s,y)\,\La(\dd s,\dd y), \quad t\in I, \]
the process $H\cdot\La=(H\cdot\La_t)_{t\in I}$ has a modification that is a semimartingale on $I$. In the case $I=\bbr$, we mean by this that $X_{-\infty} := \lim_{t\downarrow -\infty} X_t$ exists as a limit in probability, and for all bijective increasing functions $\phi\colon\bbr_+\to [-\infty,\infty)$ the process $X^\phi:=(X_{\phi(t)})_{t\in\bbr_+}$ is a usual semimartingale with respect to $(\calf_{\phi(t)})_{t\in\bbr_+}$. For later reference, we shall mention that its quadratic variation process is defined by $[X]_t:=[X^\phi]_{\phi^{-1}(t)}$ for $t\in\bar\bbr$. Finally, given a function $H\in\tilde\calp$, one can define a new random measure $H.\La$ by setting
\begin{align} K\in L^0(H.\La) &:\Leftrightarrow KH\in L^0(\La),\nonumber \\
\int_I\int_{\bbr^d} K(t,x) \,(H.\La)(\dd t,\dd x) &:= \int_I\int_{\bbr^d} K(t,x)H(t,x)\,\La(\dd t,\dd x). \label{HM} \end{align}
This indeed defines a random measure $H.\La$ if there exists a sequence $(A_k)_{k\in\bbn}\subseteq\tilde\calp$ with $A_k \uparrow \tilde\Om$ such that $\bone_{A_k}\in L^0(H.\La)$ for all $k\in\bbn$. 

Every L\'evy basis $\La$ has a canonical decomposition of the following form, see e.g. \citep[Thm.~3.2]{Chong14}: 
\beq\label{candecLa} \La(\dd t,\dd x) = B(\dd t,\dd x) + \La^\cc(\dd t,\dd x) + \int_\bbr z\bone_{\{|z|\leq1\}}\,(\mu-\nu)(\dd t,\dd x,\dd z) + \int_\bbr z\bone_{\{|z|>1\}}\,\mu(\dd t,\dd x,\dd z), \eeq
where the ingredients are as follows: 
\benu
	\item $B$ is a deterministic $\si$-finite signed Borel measure on $I\times\bbr^d$.
	\item $\La^\cc$, the continuous part of $\La$ in the usual sense (\citep[Thm.~4.13]{Bichteler83}), is a Gaussian random measure with variance measure $C$, which means that it is itself a L\'evy basis and $\La^\cc(\Om\times B)$ has a normal distribution with mean $0$ and variance $C(B)$ for every $B\in\calb_\bb$.
	\item $\mu$ is a Poisson measure on $I\times\bbr^d\times\bbr$ relative to $\bbf$ with intensity measure $\nu$, see \citep[Def.~II.1.20]{Jacod03}.
\eenu

Moreover, we have a representation
	 \begin{align} B(\dd t,\dd x)&=b(t,x)\,\la(\dd t,\dd x), &C(\dd t,\dd x)=c(t,x)\,\la(\dd t,\dd x),\nonumber \\
 \nu(\dd t, \dd x, \dd z) &= \pi(t,x,\dd z)\,\la(\dd t,\dd x), \label{charLambda} \end{align}
 with measurable functions $b\colon I\times\bbr^d\to\bbr$, $c\colon I\times\bbr^d\to\bbr_+$, a transition kernel $\pi$ from $I\times\bbr^d$ to $\bbr$ such that $\pi(t,x,\cdot)$ is a L\'evy measure for each $(t,x)$, and a positive $\si$-finite measure $\la$ on $I\times\bbr^d$ satisfying $\la(\{t\}\times\bbr^d)=0$ for all $t\in I$. 

If $\pi$ satisfies
\begin{align} \label{helppi1} \int_{|z|>1} |z|\,\pi(t,x,\dd z)&< \infty,\\
\text{or} \quad \int_{|z|\leq1} |z|\,\pi(t,x,\dd z)&< \infty,\quad\text{respectively},\label{helppi2} \end{align}
for all $(t,x)\in I\times\bbr^d$, then it makes sense to introduce the \emph{mean measure} (resp. \emph{drift measure})
\begin{align} \label{defb1} B_1(\dd t,\dd x)&:= b_1(t,x)\,\la(\dd t,\dd x), &b_1(t,x) &:= b(t,x) + \int_\bbr z\bone_{\{|z|>1\}}\,\pi(t,x,\dd z), \\
\label{defb0} B_0(\dd t,\dd x)&:=b_0(t,x)\,\la(\dd t,\dd x), &b_0(t,x) &:= b(t,x) - \int_\bbr z\bone_{\{|z|\leq1\}}\,\pi(t,x,\dd z). \end{align}
If in the first case we have $b_1(t,x)=0$ for all $(t,x)\in I\times\bbr^d$, then $\La$ is called a \emph{martingale L\'evy basis}, which will be denoted by $\La\in\calm$; if in the second case we have $b_0(t,x)=0$ for all $(t,x)\in I\times\bbr^d$, then $\La$ is called a \emph{L\'evy basis without drift}. Next, $\La$ is called \emph{symmetric} if for all $(t,x)\in I\times\bbr^d$ we have $b(t,x)=0$ and the L\'evy measure $\pi(t,x,\cdot)$ is symmetric. Furthermore, $\La$ is called a \emph{homogeneous} L\'evy basis if $\la$ is the Lebesgue measure on $I\times\bbr^d$ and $b$, $c$ and $\pi$ do not depend on $(t,x)\in I\times\bbr^d$. In this case, a function $\phi\in\tilde\calp$ is \emph{jointly stationary with $\La$} if for arbitrary $n\in\bbn$, $(h,\eta)\in \bbr\times\bbr^d$, points $(t_1,x_1),\ldots,(t_n,x_n)\in I\times\bbr^d$ and pairwise disjoint sets $B_1,\ldots,B_n\in\calb_\bb$, we have
\[ (\phi(t_i,x_i), \La(B_i) \colon i=1,\ldots, n, t_i + h \in I) \eqd (\phi(t_i+h,x_i+\eta), \La(B_i+(h,\eta)) \colon i=1,\ldots, n, t_i + h \in I). \]


Let us come back to Equation~\eqref{SPDE-var}. We first clarify what we mean by a solution $Y$ to \eqref{SPDE-var}:
\bdf\label{defsol} Equation~\eqref{SPDE-var} is said to have a \emph{solution} if there exists a predictable process $Y\in\tilde\calp$ such that for all $(t,x)\in I\times\bbr^d$ the stochastic integral on the right-hand side of \eqref{SPDE-var} is well defined and equation \eqref{SPDE-var} holds a.s. We identify two solutions $Y_1$ and $Y_2$ if for all $(t,x)\in I\times\bbr^d$ we have $Y_1(t,x)=Y_2(t,x)$ a.s.
\edf
In order to construct solutions to \eqref{SPDE-var}, we introduce some spaces of stochastic processes. Let $w\colon I\times \bbr^d \to\bbr$ be a \emph{weight function}, that is, a strictly positive measurable function. We denote by $L^{\infty,w}_I$ the Banach space of all measurable functions $f\colon I\times \bbr^d\to\bbr$ satisfying
\beq \label{Linfw} \|f\|_{L^{\infty,w}_I} := \sup_{(t,x)\in I\times \bbr^d} \frac{|f(t,x)|}{w(t,x)} <\infty.\eeq
Similarly, for $p\in(0,\infty)$, $B^{p,w}_I$ is the space of all $\phi\in\tilde\calp$ with
\beq \|\phi\|_{B^{p,w}_I} := \sup_{(t,x)\in I\times \bbr^d} \left(\frac{\bbe[|\phi(t,x)|^p]}{w(t,x)}\right)^{1/(p\vee1)} <\infty. \label{Binfw} \eeq
If $f\in L^{\infty,w}_{I_T}$ or $\phi\in B^{p,w}_{I_T}$ for all $T\in I$, then we write $f\in L^{\infty,w}_{I,\loc}$ or $\phi\in B^{p,w}_{I,\loc}$, respectively. In the special case $w\equiv1$, we use the notations $L^\infty_I$, $L^\infty_{I,\loc}$, $B^p_I$ and $B^p_{I,\loc}$.

Before we proceed to the main results of this paper, we recall how stochastic PDEs can be treated in the framework of \eqref{SPDE-var}. Let $I\subset\bbr$ be an interval, $U$ an open subset of $\bbr^d$ with boundary $\partial U$ and $P$ a polynomial in $1+d$ variables. Given some deterministic coefficient $\si$ and some L\'evy basis $\La$, they give rise to the following formal equation:
\beq\label{spde} P(\partial_t,\partial_1,\ldots,\partial_d)Y(t,x)=\si(Y(t,x))\dot\La(t,x),\quad (t,x)\in I\times U, \eeq
where $\dot\La = \partial_t\partial_1\ldots\partial_d \La$ is the formal derivative of $\La$, its noise. Usually, \eqref{spde} is subjected to some boundary conditions on $\partial (I\times U)$. Of course, the derivative of $\La$ is not defined except in trivial cases, so a strong solution to \eqref{spde} will not exist. Going back to \citep{Walsh86} is the idea of constructing a so-called \emph{mild solution} to \eqref{spde}. For this method to work, one has to assume that the operator $P$ possesses a Green's function on $I\times U$. Then a mild solution to \eqref{spde} is nothing but a solution in the sense of Definition~\ref{defsol} to \eqref{SPDE-var}, where $G$ is the Green's function and $Y_0$ a term that only depends on the boundary conditions posed on $\partial (I\times U)$.

\brem
While the notion of a solution as in Definition~\ref{defsol} is very common in the theory of stochastic PDEs, it is different to the standard notion of solutions to (ordinary) SDEs: let $I=\bbr_+$ and $d=0$, that is, space contains only one point, and consider $G(t,1;s,1)=g(s)\bone_{\{s\leq t\}}$ with some smooth function $g$. Then Equation~\eqref{SPDE-var} is equivalent to the SDE
\beq\label{SDEex} \dd Y(t) = g(t)\si(Y(t-))\,\La(\dd t), \quad t\geq0, \quad Y(0)=Y_0, \eeq
where $\La$ is a semimartingale with independent increments. Ordinary SDE theory tells us that Equation~\eqref{SDEex} has a c\`adl\`ag solution $Y$ that is unique up to indistinguishability. In contrast, a solution in the sense of Definition~\ref{defsol} would be the predictable version $Y(\cdot-)$, and uniqueness is only understood up to modifications. The reason why we have chosen this slightly different notion of a solution is that we are particularly interested in the case where $G$ in Equation~\eqref{SPDE-var} has singularities. In such cases, Equation~\eqref{SPDE-var} permits no c\`adl\`ag solutions.
\erem

\section{Existence and uniqueness results on $I=\bbr_+$}\label{localex}

The goal of this section is to provide sufficient conditions under which there exists a (unique) solution to \eqref{SPDE-var} on the interval $I=\bbr_+$. It is clear that everything in this section holds analogously if we replace $I=[0,\infty)$ by $I=[a,\infty)$ with some $a\in\bbr$. As mentioned in the Introduction, the forthcoming theorem generalizes the results of \citep{SLB98} to potentially inhomogeneous L\'evy bases and kernels different from the heat kernel. It holds under the following list of assumptions:

\bass\label{Aass} Let $p\in(0,2]$ and the predictable characteristics of $\La$ be given by \eqref{charLambda}. We impose the following conditions:
\benu
\item $Y_0\in B^p_{[0,\infty),\loc}$.
\item There exists $C_{\si,1}\in\bbr_+$ such that $|\si(x)-\si(y)|\leq C_{\si,1} |x-y|$ for all $x, y \in \bbr$.
\item $G\colon (\bbr_+\times\bbr^d)^2 \to \bbr$ is a measurable function such that $G(t,\cdot;s,\cdot)\equiv0$ whenever $s>t$.
\item If $p<2$, then $\La$ has no Gaussian part: $c(t,x)=0$ for $(t,x)\in\bbr_+\times\bbr^d$. If $p=2$, then we assume for all $T\in\bbr_+$
\beq\label{Gaussianfinite} \sup_{(t,x)\in[0,T]\times\bbr^d} \int_0^t\int_{\bbr^d} |G(t,x;s,y)|^2 c(s,y)\,\la(\dd s,\dd y) < \infty. \eeq
\item For all $T\in\bbr_+$
\beq\label{jumpsfinite} \sup_{(t,x)\in [0,T]\times\bbr^d} \int_0^t \int_{\bbr^d} \int_\bbr |G(t,x;s,y)z|^p \,\nu(\dd s,\dd y,\dd z) < \infty. \eeq
\item Recall the definition of $b_1$ and $b_0$ from \eqref{defb1} and \eqref{defb0}. If $p\geq1$, assume that $\nu$ satisfies \eqref{helppi1} and that for all $T\in\bbr$
	\beq\label{meanstable0} \sup_{(t,x)\in [0,T]\times\bbr^d} \int_0^t \int_{\bbr^d} |G(t,x;s,y)b_1(s,y)| \,\la(\dd s,\dd y) < \infty; \eeq
if $p<1$, assume that $\nu$ satisfies \eqref{helppi2} and that $b_0(t,x)=0$ for all $(t,x)\in \bbr_+\times\bbr^d$.
\item Define for $(t,x),(s,y)\in \bbr_+\times\bbr^d$
\[ G^{A}(t,x;s,y):= |G(t,x;s,y)|^p\left(\int_\bbr |z|^p \,\pi(s,y,\dd z) + c(s,y)\right) + |G(t,x;s,y)b_1(s,y)|\bone_{\{p\geq1\}}, \]
and assume that for every $T\in \bbr_+$ and $\eps>0$ there exists $k\in\bbn$ together with a subdivision $\mathcal{T}\colon 0=t_0 < t_1 < \ldots < t_{k+1} =T$ such that
\beq\label{partG0} \sup_{(t,x)\in [0,T]\times\bbr^d} \sup_{i=0,\ldots,k} \int_{t_i}^{t_{i+1}}\int_{\bbr^d} G^{A}(t,x;s,y)\,\la(\dd s,\dd y)<\eps. \eeq
\eenu
\eass

\bthm\label{main1} 
Let Assumption~\ref{Aass} be valid. Then Equation \eqref{SPDE-var} has a unique solution in $B^p_{[0,\infty),\loc}$. 
\ethm

The conditions of Assumption~\ref{Aass} simplify a lot if $G$ and $\La$ are \emph{quasi-stationary}, that is,
\beq\label{quasistat}
	|G(t,x;s,y)|\leq g(t-s,x-y),\quad \la(\dd t,\dd x)=\dd(t,x),\quad b, c \in L^\infty_{[0,\infty),\loc},\quad \pi(t,x,\dd z) \leq \pi_0(\dd z), 
\eeq 
where $g\colon \bbr_+\times\bbr^d\to\bbr$ is a positive measurable function. 
\bcor\label{cor1} Suppose that \eqref{quasistat} holds and that Assumption~\ref{Aass}(1), (2) and (3) are given. Furthermore, assume that we have for some $p\in(0,2]$
\beq\label{reg} b_0\equiv0 \text{ if } p<1,\quad c\equiv0 \text{ if } p<2,\quad \int_\bbr |z|^p\,\pi_0(\dd z)<\infty, \eeq
and for all $T\in\bbr_+$
\beq\label{localint}
\int_0^T\int_{\bbr^d} g^p(t,x)\,\dd(t,x) < \infty,\quad\text{and}\quad \int_0^T\int_{\bbr^d} g(t,x)\,\dd(t,x) < \infty \text{ if } p\geq1 \text{ and } \La\notin\calm.
\eeq
Then all conditions of Assumption~\ref{Aass} are satisfied and Theorem~\ref{main1} holds.
\ecor

\brem\label{rem1}
\benu
	\item Assumption~\ref{Aass} and Theorem~\ref{main1} are special cases of Assumption~\ref{Cass} and Theorem~\ref{main3}, respectively, which we will discuss in Section~\ref{globalex}. In fact, Theorem~\ref{main1} follows if we take $I=[0,\infty)$ and $w\equiv1$ in Theorem~\ref{main3}.
	\item Conditions (4), (5) and (6) in Assumption~\ref{Aass} are conditions on the joint size of $G$ and the three characteristics of $\La$, respectively. Although they are valid for many interesting examples, especially condition (5) might be too restrictive: it is violated as soon as the moment structure of $\La$ is worse than its variation structure, which, for instance, occurs if $\La$ is an $\al$-stable L\'evy basis with $\al\in(0,2)$; see also the last condition in \eqref{reg}. Theorem~\ref{main2} below provides, under some additional hypotheses, an extension of Theorem~\ref{main1} that includes such cases.
	\item The following observation follows from Corollary~\ref{cor1}: in the quasi-stationary case \eqref{quasistat}, condition (7) in Assumption~\ref{Aass} is already implied by conditions (4), (5) and (6). In other words, condition (7) is a smallness assumption on the non-stationary part of $G$ and the characteristics of $\La$.
	\item As we shall see in the more general Theorem~\ref{main3} in Section~\ref{globalex}, it actually suffices that the left-hand side of \eqref{partG0} can be made smaller than some fixed constant that does not depend on $\mathcal{T}$. Due to the previous remark, however, this fact is not that important in the case $I=[0,\infty)$ (in the case $I=\bbr$, it is!).
	\eenu
\erem

Next, we apply Theorem~\ref{main1} and its corollary to the stochastic heat equation. In fact, this equation will serve as our toy example and will be carried through the whole paper and revisited after each main theorem: see the Examples~\ref{ex2}, \ref{ex3} and \ref{ex4}.
\bex\label{ex1} We consider the stochastic heat equation on $\bbr_+\times\bbr^d$, that is, \eqref{spde} with $P$ given by $P(t,x)=t-\sum_{i=1}^d x_i^2 + a$, $a\in\bbr$, and some Lipschitz coefficient $\si$. The Green's function is the heat kernel
\beq\label{heatkernel} G_a(t,x;s,y)=g_a(t-s,x-y)=\frac{\exp\left(-\frac{|x-y|^2}{4(t-s)}-a(t-s)\right)}{(4\pi(t-s))^{d/2}}\bone_{\{s<t\}}. \eeq
We pose an initial condition at time $t=0$, that is, we require $Y(0,x)=y_0(x)$, where $y_0\colon\bbr^d\to\bbr$ is some bounded continuous and, for simplicity, deterministic function. Then the correct term for $Y_0$ in \eqref{SPDE-var} is
\beq\label{initial} Y_0(t,x):=\int_{\bbr^d} g_a(t,x-y)y_0(y)\,\dd y,\quad (t,x)\in\bbr_+\times\bbr^d. \eeq
The stochastic heat equation on $I=\bbr_+$ then reads as
\beq\label{heat1} Y(t,x)=Y_0(t,x)+ \int_0^t\int_{\bbr^d} g_a(t-s,x-y)\si(Y(s,y))\,\La(\dd s,\dd y),\quad (t,x)\in\bbr_+\times\bbr^d.\eeq

Let us determine sufficient conditions for existence and uniqueness of solutions to \eqref{heat1}: assuming that the characteristics of $\La$ satisfy \eqref{quasistat}, we have to check the conditions of Corollary~\ref{cor1}: (1) and (2) of Assumption~\ref{Aass} are clear. Since 
\beq\label{intheatloc} \int_0^T \int_{\bbr^d} g_a^p(s,y)\,\dd(s,y) < \infty \text{ for all } T\in\bbr_+ \iff p<1+2/d,\eeq
we obtain existence and uniqueness for the stochastic heat equation \eqref{heat1} on $I=\bbr_+$ if \eqref{reg} holds with some $0<p<1+2/d$. In particular, this excludes the choice $p=2$ and therefore the possibility of taking a non-zero Gaussian part whenever $d\geq2$.
\eex

As pointed out in Remark~\ref{rem1}(2), Theorem~\ref{main1} excludes any L\'evy basis that has the property that for every $p\in(0,2]$ 
\beq\label{badmom} \la\left(\left\{(t,x)\in\bbr_+\times\bbr^d\colon \int_\bbr |z|^p\,\pi(t,x,\dd z) = \infty\right\}\right)>0. \eeq
We now discuss a possibility to circumvent this. 

\bass\label{Bass} Consider the following hypotheses:
\benu
	\item Assumption~\ref{Aass}(2) and (3) are valid. 
	\item There exists some $q\in(0,2]$ such that for all $n\in\bbn$ conditions (4)--(7) of Assumption~\ref{Aass} are valid when $p$ is replaced by $q$ and $\nu$ is replaced by 
	\[ \nu^n(\dd t,\dd x,\dd z) := \bone_{\{|z|\leq n\}}\,\nu(\dd t,\dd x,\dd z). \]
	Of course, $b_1$ is changed accordingly.
	\item For all $T\in\bbr_+$ we have $\nu\big([0,T]\times\bbr^d\times[-1,1]^\Comp\big)<\infty$.
	\item $Y_0\in \tilde\calp$ and there are stopping times $(T_n)_{n\in\bbn}$ with $T_n\uparrow\infty$ a.s. and $Y_0\bone_{\llbracket 0,T_n\rrbracket}\in B^q_{[0,\infty),\loc}$ for all $n\in\bbn$.	
	\item There exist $\ga\in(0,1)$ and $C_{\si,2}\in\bbr_+$ such that $|\si(x)|\leq |\si(0)|+C_{\si,2}|x|^\ga$ for all $x\in\bbr$.
	\item There exists $p\in(0,2)$ satisfying $p<q$ and $q\ga\leq p$ such that $Y_0\in B^p_{[0,\infty),\loc}$.
	\item For all $T\in\bbr_+$
	\[ \sup_{(t,x)\in[0,T]\times\bbr^d} \int_0^t \int_{\bbr^d} \int_\bbr |G(t,x;s,y)z|^p_q \,\nu(\dd s,\dd y,\dd z) < \infty. \]
	\item If $p\geq1$, \eqref{meanstable0} holds.
	\item If $p<1$, there exist exponents $\al\in(-\infty,2],\beta\in[0,\infty)$ with the following properties:
	\item[(9a)] For all $(t,x)\in\bbr_+\times\bbr^d$, $A\in[1,\infty)$ and $a\in(0,1]$ we have
	\begin{align}	\left|b(t,x) - \int_\bbr z\bone_{\{|z|\in(a,1]\}} \,\pi(t,x,\dd z)\right| &\leq F_0(t,x) a^{1-\al}, \label{driftconv}\\
	 \left| b(t,x) + \int_\bbr z\bone_{\{|z|\in(1,A]\}} \,\pi(t,x,\dd z)\right| &\leq F_1(t,x) A^{1-\beta} \label{meanconv}
	\end{align}
	for some positive measurable functions $F_0,F_1\colon \bbr_+\times\bbr^d\to\bbr$.
	\item[(9b)] For all $T\in\bbr_+$ we have
	\beq\label{Ginteg} \sup_{(t,x)\in[0,T]\times\bbr^d} \int_0^t\int_{\bbr^d} (F_0(s,y)\vee F_1(s,y))|G(t,x;s,y)|^\al_\beta\,\la(\dd s,\dd y) < \infty. \eeq
	\item[(9c)] $(\al\vee\beta)\ga\leq p$.
	\item[(10)] The partition property \eqref{partG0} holds with $G^{B}$ instead of $G^{A}$, where for $(t,x), (s,y)\in\bbr_+\times\bbr^d$
	\begin{align} G^{B}(t,x;s,y)&:=|G(t,x;s,y)|^2 c(s,y) + \int_\bbr |G(t,x;s,y)z|^p_q\,\pi(s,y,\dd z)  \nonumber\\
&\quad+ \begin{cases} |G(t,x;s,y)b_1(s,y)|, &p\geq1,\\ (F_0(s,y)\vee F_1(s,y))|G(t,x;s,y)|^\al_\beta, &p<1 \end{cases}, \label{GB} \end{align}
	\eenu
\eass

\bthm\label{main2} \benu \item Suppose that conditions (1)--(4) of Assumption~\ref{Bass} are true. Then there exists a unique solution to Equation~\eqref{SPDE-var} among those $Y\in \tilde\calp$ for which there exist stopping times $(T_n)_{n\in\bbn}$ with $T_n\uparrow\infty$ a.s. such that $Y\bone_{\llbracket 0,T_n\rrbracket} \in B^q_{[0,\infty),\loc}$ for all $n\in\bbn$. \item If in addition also conditions (5)--(10) of Assumption~\ref{Bass} are valid, then the solution $Y$ from part (1) belongs to $B^p_{[0,\infty),\loc}$. \eenu
\ethm

\brem 
\benu 
	\item Part (1) of this theorem relies on some stopping time techniques that have already been used in \citep{Balan14} to construct solutions to \eqref{SPDE-var} driven by $\al$-stable noise with $\al\neq1$. Theorem~\ref{main2} extends this result to more general L\'evy bases and, more importantly, provides in part (2) conditions under which this solution belongs to the space $B^p_{[0,\infty),\loc}$. 
	\item The smaller the growth index $\ga$ of $\si$ is, the smaller can $p$ be chosen and therefore, the weaker the conditions (6)--(9) of Assumption~\ref{Bass} are. For $\al$-stable L\'evy bases with $\al\in(0,2)$, any $\ga\in(0,1)$ and $p<q$ will suffice.
	\item If $p<1$, condition (9) of Assumption~\ref{Bass} looks quite technical but is actually only a very mild assumption. In the next Corollary~\ref{cor2} where we treat the quasi-stationary case, it is already implied by condition (6) below.
	\item Remark~\ref{rem1}(3) holds analogously: see the next corollary.
	\item For the second condition of Assumption~\ref{Bass}, if $p\geq1$, one has to check Assumption~\ref{Aass}(6) for different replacements of $b_1$ as $n$ varies, which are usually non-zero even when $\La\in\calm$.
	\item The most stringent condition in Assumption~\ref{Bass} is (3): it requires the intensity of large jumps of $\La$ to decay quickly enough in space. For example, it is typically \emph{not} enough to have $\pi(t,x,\dd z)=\pi_0(\dd z)$. See Corollary~\ref{cor2} and Example~\ref{ex2} for more details.
\eenu
\erem

Again we reformulate Assumption~\ref{Bass} in the quasi-stationary case:
\bcor\label{cor2} Assume that $G$ and $\La$ satisfy \eqref{quasistat}, but with the stronger condition
\beq\label{spacedecay} \pi(t,x,\dd z) \leq \pi_1(t,x)\,\pi_0(\dd z),\quad \pi_1\in L^\infty_{[0,\infty),\loc},\quad \int_0^T\int_{\bbr^d} \pi_1(t,x)\,\dd(t,x)<\infty \eeq
for all $T\in\bbr_+$. Then part (1) of Theorem~\ref{main2} holds if:
\benu
	\item Assumption~\ref{Bass}(1) and (4) are valid.
	\item For some $q\in(0,2]$ conditions \eqref{reg} and \eqref{localint} hold with $p$ replaced by $q$ and $\pi_0$ replaced by $\bone_{\{|z|\leq1\}}\,\pi_0(\dd z)$.
	\item If $q\geq1$, either $\int_0^T\int_{\bbr^d} g(t,x)\,\dd(t,x)<\infty$ for all $T\in\bbr_+$, or $\La$ is symmetric.
\eenu
Part (2) of the same theorem holds if additionally:
\benu
	\item[(4)] $\si$ satisfies the growth condition of Assumption~\ref{Bass}(5) with $\ga\in(0,1)$. 
	\item[(5)] There exists $p\in(0,2)$ with $p<q$ and $q\ga\leq p$ such that $Y_0\in B^p_{[0,\infty),\loc}$.
	\item[(6)] $\displaystyle \int_\bbr |z|^p_q\,\pi_0(\dd z) < \infty$ and $\displaystyle \int_0^T\int_{\bbr^d} |g(t,x)|^q_p\,\dd(t,x) < \infty$ for all $T\in\bbr_+$.
\eenu
\ecor 

For illustration purposes we go through the conditions of Theorem~\ref{main2} and Corollary~\ref{cor2} for the stochastic heat equation.
\bex[Continuation of Example~\ref{ex1}]\label{ex2}
Our aim is to extend the findings of Example~\ref{ex1} when $\La$ has bad moment properties in the sense of \eqref{badmom}. For simplicity we assume that the characteristics of $\La$ are within the setting of Corollary~\ref{cor2}, that is, they satisfy \eqref{quasistat} and \eqref{spacedecay}. As before, $\si$ is a Lipschitz continuous function and the equation of interest is \eqref{heat1} with $Y_0$ given by \eqref{initial}. In view of \eqref{intheatloc}, it is immediate to see that Corollary~\ref{cor2} yields the following conditions for part (1) of Theorem~\ref{main2} to hold:
\beq\label{B1-cor} \int_{[-1,1]} |z|^q\,\pi_0(\dd z) <\infty \text{ for some } 0<q<1+2/d,\quad c\equiv0 \text{ if } d\geq2,\quad b_0\equiv0 \text{ if } q<1. \eeq

Furthermore, if $\si$ has growth of order $\ga\in(0,1)$ and
\beq\label{B2-cor} \int_{|z|>1} |z|^p\,\pi_0(\dd z) < \infty \text{ for some } p<1+2/d \text{ with } p<q \text{ and } q\ga\leq p, \eeq
then the solution $Y$ belongs to $B^p_{[0,\infty),\loc}$. Indeed, this claim follows from Corollary~\ref{cor2} and the fact that for all $p,q\in(0,\infty)$ we have
\beq\label{pgG} \int_0^T\int_{\bbr^d} |g_a(t,x)|^q_p \,\dd(t,x) < \infty \eeq
for all $T\in\bbr_+$ if and only if $q\in(0,1+2/d)$ ($p$ does not matter). From \eqref{B2-cor} we also see the following: the smaller the growth order $\ga$ of $\si$ is, the fewer moments $\pi_0$ is required to have.

At last, we give some further explanation for the integrability condition on $\pi_1$ given in \eqref{spacedecay}. We assume that $\pi(t,x,\dd z)=\pi_1(t,x)\pi_0(\dd z)$ with a L\'evy measure $\pi_0$ of unbounded support. Then it is obvious to see that we cannot take $\pi_1\equiv1$, that is, a homogeneous noise $\La$, but have to choose $\pi_1$ with sufficient decay in space. For instance, if there exists some exponent $r\in\bbr$ such that for all $T\in\bbr_+$ we have $\pi_1(t,x)\leq C_T|x|^{-r}$ for all $(t,x)\in [0,T]\times\bbr^d$ and some constant $C_T\in\bbr_+$, then we need for \eqref{spacedecay} that $r>d$, a condition that is stronger in higher dimensions. Finally, \eqref{spacedecay} is always met if $\pi_1$ is bounded and vanishes outside a compact in $\bbr^d$, which corresponds to a noise that only acts locally. In particular, this assumption is very natural if we consider the stochastic heat equation on bounded domains as, for instance, in \citep{Applebaum00,Balan14,Walsh86}.
\eex

\brem Theorem~\ref{main1} and \ref{main2} can actually be extended to even more general random measures than L\'evy bases. Let us consider a random measure $M$ on $\bbr_+\times\bbr^d$ that is defined by
\begin{align} M(\dd t,\dd x) &= b(t,x)\,\dd(t,x) + \rho(t,x)\,W(\dd t,\dd x) + \int_E \un\delta(t,x,z)\,(\pf-\qf)(\dd t,\dd x,\dd z) \nonumber\\
&\quad + \int_E \ov\delta(t,x,z)\,\pf(\dd t,\dd x,\dd z), \label{candecM} \end{align}
where $(E,\cale)$ is an arbitrary Polish space equipped with its Borel $\si$-field, $b,\rho\in\tilde\calp$, $\delta=\un\delta+\ov\delta=\delta\bone_{\{|\delta|\leq 1\}} + \delta\bone_{\{|\delta|>1\}}$ is an $\tilde\calp\otimes\cale$-measurable function, $W$ is a Gaussian random measure with the Lebesgue measure on $I\times\bbr^d$ as variance measure, $\pf$ is a homogeneous Poisson random measure on $I\times\bbr^d\times E$ relative to the filtration $\bbf$ with intensity measure $\qf(\dd t,\dd x,\dd z) = \dd t\,\dd x\,\la(\dd z)$ where $\la$ is a $\si$-finite infinite atomless measure on $(E,\cale)$. Moreover, all ingredients are such that $M(\Om\times (I\cap(-k,k])\times (-k,k]^d)$ is well defined for all $k\in\bbn$. Such a measure $M$ can be viewed as the space--time analogue of It\^o semimartingales. We impose the following conditions on the coefficients (these are classical in the semimartingale setting, cf. \citep[Chap.~6]{AitSahalia14}): there exist positive constants $(\beta_N)_{N\in\bbn}$, a sequence of stopping times $(\tau_N)_{N\in\bbn}$ increasing to infinity a.s., and deterministic positive measurable functions $j_N(z)$ such that for all $(\om,t,x)\in\tilde\calp$ with $t\leq \tau_N(\om)$ we have
	\benu
		\item $|b(\om,t,x)|, |c(\om,t,x)|\leq \beta_N$, 
		\item $|\delta(\om,t,x,z)|^p \leq j_N(z)$ and $\int_E j_N(z)\,\la(\dd z)<\infty$.
	\eenu
Then with obvious changes to Assumptions~\ref{Aass} and \ref{Bass}, respectively, Theorems~\ref{main1} and \ref{main2} also apply to Equation~\eqref{SPDE-var} when driven by the random measure $M$ as given in \eqref{candecM}.
\erem

\section{Existence and uniqueness results on $I=\bbr$}\label{globalex}

While Section~\ref{localex} deals with Equation~\eqref{SPDE-var} on $I=[0,\infty)$, this section investigates the case $I=\bbr$. In particular, we obtain conditions for Equation~\eqref{SPDE-var} to possess a stationary solution. In order to demonstrate the difference between the two cases $I=[0,\infty)$ and $I=\bbr$, we analyse the following deterministic example.

\bex\label{exexpfunc} Let $\la\in\bbr$ and consider the following equation:
\beq\label{expex} v(t) =  1+\int_{-\infty}^t \ee^{-\la(t-s)}v(s)\,\dd s,\quad t\in\bbr. \eeq
By standard computation one can show the following: if $\la\leq0$, Equation \eqref{expex} has no solution; if $\la>0$ and $\la\neq1$, then the solutions to \eqref{expex} are 
\[ v(t)= c\ee^{(1-\la)t}+\frac{\la}{\la-1},\quad c\in\bbr; \]
if $\la=1$, the solutions are
\[ v(t) = t + c,\quad c\in\bbr.\]
We draw some important conclusions, also regarding possibilities and limitations for Equation~\eqref{SPDE-var} with $I=\bbr$:
\benu
	\item The reason why \eqref{expex} possesses no solution for $\la\leq0$ is simply the non-integrability of the kernel:
	\beq\label{kernelint} \int_{-\infty}^t \ee^{-\la(t-s)}\,\dd s = \int_0^\infty \ee^{-\la s} \,\dd s = \infty. \eeq
	\item If Equation~\eqref{expex} has a solution, it has uncountably many. If $\la\in(1,\infty)$, only one solution is in $L^{\infty}_{\bbr,\loc}$, namely if $c=0$. The reason for this is that the integral of the kernel given in \eqref{kernelint} is smaller than $1$. In this case the uniqueness of solutions in $L^{\infty}_{\bbr,\loc}$ follows from Lemma~\ref{comparison}(2). Thus, in the stochastic case of \eqref{SPDE-var}, we can expect existence and uniqueness of solutions in $B^p_{\bbr,\loc}$ only if the quantities \eqref{Gaussianfinite}, \eqref{jumpsfinite} and \eqref{meanstable0} are \emph{small enough} (not only finite) in a sense to be made precise.
	\item In contrast to the case $\la\in(1,\infty)$, we have for $\la\in(0,1)$ that all solutions belong to $L^\infty_{\bbr,\loc}$ and for $\la=1$ that no solution belongs to $L^\infty_{\bbr,\loc}$. Furthermore, in these cases, all solutions start with strictly negative values at $-\infty$. This is somewhat surprising given the fact that all ingredients of \eqref{expex} (the exponential kernel, the constant driving force and the Lebesgue measure as integrator) are positive. This phenomenon is typical when the integral of the kernel in \eqref{kernelint} becomes greater or equal to one: the kernel is too large to allow for a positive solution. Finally, none of the solutions can be found via a Picard iteration scheme (since the Picard iterates are always positive when the input factors are). Thus, if the kernel in \eqref{SPDE-var} is too large in a certain sense, we will not be able to construct a solution in general.
	\item Under certain circumstances, however, one can make the kernel size smaller (which then implies the existence and uniqueness of solutions) by considering Volterra equations in weighted spaces. For instance, consider the following modification of Equation~\eqref{expex}:
	\beq\label{expex2} v(t) = \ee^{\al t} + \int_{-\infty}^t \ee^{-\la(t-s)}v(s)\,\dd s,\quad t\in\bbr,\eeq
	with $\al,\la\in\bbr$ satisfying $\la>0$ and $\al+\la>1$. The family of solutions in this case is
	\beq\label{sol2} v(t) = \frac{\al+\la}{\al+\la-1}\ee^{\al t} + c\ee^{(1-\la)t},\quad t\in\bbr,\quad c\in\bbr. \eeq
	First note that positive solutions do exist, namely, when $c\geq0$. Furthermore, with $w(t):= \ee^{\al t}$, we have
	\[ \int_{-\infty}^t w^{-1}(t)\ee^{-\la(t-s)}w(s)\,\dd s= \int_{-\infty}^t \ee^{-(\al + \la)(t-s)}\,\dd s = (\al+\la)^{-1}<1.  \]
	That is, by Lemma~\ref{comparison}(2), there exists a unique solution to \eqref{expex2} in $L^{\infty,w}_{\bbr,\loc}$, which corresponds to the case $c=0$ in \eqref{sol2}. Roughly speaking, this device was possible because the force function $\ee^{\al t}$ is small enough at $-\infty$ (the constant function in \eqref{expex} was obviously \emph{not} small enough). This motivates us to work in the weighted spaces $B^{p,w}_{\bbr,\loc}$ for Equation~\eqref{SPDE-var} on $I=\bbr$.
\eenu
\eex

We are about to formulate a set of conditions that generalizes those of Assumption~\ref{Aass} and leads to the existence and uniqueness of solutions for Equation~\eqref{SPDE-var} on arbitrary intervals, in particular on $I=\bbr$. In order to do so, we need the following definition.
\bdf\label{BDGineq} Let $p\in(0,\infty)$. \benu \item For $p\in(0,1)$ we set $C^\mathrm{BDG}_p:=1$. \item For $p\in[1,\infty)$ we denote by $C^\mathrm{BDG}_p$ the smallest positive number such that for all local martingales $(M_t)_{t\in\bbr_+}$ w.r.t. $\bbf$ we have 
\beq\label{BDG} \sup_{t\geq0} \|M_t\|_{L^p} \leq C^{\mathrm{BDG}}_p \|[M]^{1/2}_\infty\|_{L^p}. \eeq \eenu
\edf

\brem We make some comments on Definition~\ref{BDGineq}: \benu
	\item The Burkholder-Davis-Gundy inequality ensures the finiteness of $C^\mathrm{BDG}_p$ for $p\in[1,\infty)$. Of course, inequality \eqref{BDG} becomes false in general for $p<1$; the definition above for $p\in(0,1)$ is merely for notational convenience. Moreover, the inequality for $p\in[1,\infty)$ is usually stated with the supremum inside the $L^p$-norm on the left-hand side of \eqref{BDG}. However, this may enlarge the optimal constant $C^\mathrm{BDG}_p$.
	\item The choice $I=\bbr_+$ is unimportant: a straightforward time change argument shows that $C^\mathrm{BDG}_p$ remains optimal for any other non-trivial interval $I\subseteq\bbr$.
	\item For $p\in[1,\infty)$, the actual value of $C^\mathrm{BDG}_p$ is not known in general. We are only interested in the case $p\in[1,2]$, for which the following results are available: $C^\mathrm{BDG}_p\leq \sqrt{8p}$ for $p\in(1,2)$, $C^\mathrm{BDG}_2 = 1$ (cf. \citep[Eq.~(4.2.3)]{Bichteler02}) and $C^\mathrm{BDG}_1=2$ (cf. \citep[Thm.~8.7]{Osekowski12}).\eenu
\erem

\bass\label{Cass} Let $0< p \leq 2$, $I\subseteq\bbr$ be an interval and $w\colon I\times\bbr^d\to\bbr$ be a weight function. We impose the following conditions:
\benu
\item $Y_0\in B^{p,w}_{I,\loc}$.
\item There exists $C_{\si,1}\in\bbr_+$ such that $|\si(x)-\si(y)|\leq C_{\si,1} |x-y|$ for all $x, y \in \bbr$.
\item $G\colon (I\times\bbr^d)^2 \to \bbr$ is a measurable function such that $G(t,\cdot;s,\cdot)\equiv0$ whenever $s>t$.
\item If $p<2$, then $\La$ has no Gaussian part: $c(t,x)=0$ for all $(t,x)\in I\times\bbr^d$. If $p=2$, then we assume for all $T\in I$
\beq\label{Gaussianfin} \sup_{(t,x)\in I_T\times\bbr^d} \int_I\int_{\bbr^d} w^{-1}(t,x)|G(t,x;s,y)|^2c(s,y)(w(s,y)\vee\si(0))\,\la(\dd s,\dd y) < \infty. \eeq
\item For all $T\in I$
\beq\label{nustable} \sup_{(t,x)\in I_T\times\bbr^d} \int_I \int_{\bbr^d} \int_\bbr w^{-1}(t,x)|G(t,x;s,y)z|^p (w(s,y)\vee \si(0)) \,\nu(\dd s,\dd y,\dd z) < \infty. \eeq
\item Recall the definition of $b_1$ and $b_0$ from \eqref{defb1} and \eqref{defb0}. If $p\geq1$, assume that $\nu$ satisfies \eqref{helppi1} and that for all $T\in I$
	\beq\label{meanstable} \sup_{(t,x)\in I_T\times\bbr^d} \int_I \int_{\bbr^d} w^{-1}(t,x)|G(t,x;s,y)b_1(s,y)| (w(s,y)\vee\si(0))\,\la(\dd s,\dd y) < \infty; \eeq
if $p<1$, assume that $\nu$ satisfies \eqref{helppi2} and that $b_0(t,x)=0$ for all $(t,x)\in I\times\bbr^d$.
\item If $p\geq1$ and $\La\notin\calm$, assume that (6) also holds with $w$ replaced by the constant function $1$.
\item Define for $(t,x),(s,y)\in I\times\bbr^d$
\begin{align} \bar G^{C,1}(t,x;s,y)&:= (C_{\si,1}C^\mathrm{BDG}_p)^p |G(t,x;s,y)|^p \left(\int_\bbr |z|^p \,\pi(s,y,\dd z) + c(s,y)\right), \nonumber\\
\bar G^{C,2}(t,x;s,y)&:= C_{\si,1}^p \left(\int_I\int_{\bbr^d} |G(t,x;s,y)b_1(s,y)|\,\la(\dd s,\dd y)\right)^{p-1}|G(t,x;s,y)b_1(s,y)|\bone_{\{p\geq1\}}, \nonumber\\
G^{C,1}(t,x;s,y) &:= w^{-1}(t,x)\bar G^{C,1}(t,x;s,y)w(s,y),\nonumber\\
G^{C,2}(t,x;s,y) &:= w^{-1}(t,x)\bar G^{C,2}(t,x;s,y)w(s,y),\label{deftildeG} \end{align}
and assume that for every $T\in I$ there exists $k\in\bbn$ together with a subdivision $\mathcal{T}\colon \inf I=t_0 < t_1 < \ldots < t_{k+1} =T$ such that
\beq\label{partG} \sup_{(t,x)\in I_T\times\bbr^d} \sup_{i=0,\ldots,k} \sum_{l=1}^2 \left(\int_{t_i}^{t_{i+1}}\int_{\bbr^d}  G^{C,l}(t,x;s,y)\,\la(\dd s,\dd y)\right)^{1/(p\vee1)}<1. \eeq
\eenu
\eass

\bthm\label{main3} Under Assumption~\ref{Cass} there exists a unique solution to Equation~\eqref{SPDE-var} in $B^{p,w}_{I,\loc}$. \ethm

In the quasi-stationary case, Assumption~\ref{Cass} simplifies a lot:

\bcor\label{cor3}
Let $I=\bbr$, $w\equiv1$ and Assumption~\ref{Cass}(1), (2) and (3) be valid. We assume that $G$ and $\La$ satisfy 
\beq
	|G(t,x;s,y)|\leq g(t-s,x-y),\quad \la(\dd t,\dd x)=\dd(t,x),\quad b,c\in L^\infty_\bbr,\quad \pi(t,x,\dd z) \leq  \pi_0(\dd z)
 \label{quasistat2}
\eeq
for all $(t,x)\in \bbr\times\bbr^d$ and some positive measurable $g\colon \bbr_+\times\bbr^d\to\bbr$. Furthermore, we suppose that for some $p\in(0,2]$ we have 
\beq\label{reg2} b_0\equiv0\text{ if } p<1, \quad c\equiv0 \text{ if } p<2,\quad \zeta_p:=\int_\bbr |z|^p\,\pi_0(\dd z)< \infty, \eeq
and that the following size condition is fulfilled: if $p\in(0,1)$, then
\beq\label{intconst1} C_{\si,1}^p\zeta_p\int_0^\infty \int_{\bbr^d} g^p(t,x)\,\dd(t,x)<1,\eeq
and if $p\in[1,2]$, then
\beq\label{intconst2} C_{\si,1}\Bigg[C^\mathrm{BDG}_p\left((\zeta_p + \|c\|_{L^\infty_\bbr}) \int_0^\infty\int_{\bbr^d} g^p(t,x)\,\dd(t,x)\right)^{1/p} + \|b_1\|_{L^\infty_\bbr} \int_0^\infty\int_{\bbr^d} g(t,x)\,\dd(t,x) \Bigg]<1.\eeq
Then all conditions of Assumption~\ref{Cass} are satisfied and Theorem~\ref{main3} holds.
\ecor

We write down some important observations:
\brem
\benu
	\item There is a fundamental difference between condition (7) of Assumption~\ref{Aass} and condition (8) of Assumption~\ref{Cass}. For instance, consider the quasi-stationary case in Corollary~\ref{cor1} and Corollary~\ref{cor3}, where they reduce to \eqref{localint} and either \eqref{intconst1} or \eqref{intconst2}. While in the former case we only need certain integrability properties of the kernel, we explicitly have to care about the size of the integrals involved in the latter case, which is also the size condition we have mentioned in Example~\ref{exexpfunc}(2). Also notice that this is related to the fact that in the case $I=\bbr$, we typically cannot make the left-hand side of \eqref{partG} as small as we want by refining the subdivision $\mathcal{T}$ since the first interval $(t_0,t_1]=(-\infty,t_1]$ always has infinite length. So whereas condition (7) of Assumption~\ref{Aass} is quite natural for $I=[0,\infty)$, the analogous condition for $I=\bbr$ would be very restrictive.
	\item By the nature of Equation~\eqref{SPDE-var}, the size condition (8) of Assumption~\ref{Cass} is ``symmetric'' in $G$, $\si$ and $\La$.
	\item In Theorem~\ref{main3} uniqueness does not hold in $\tilde\calp$: see Equation~\eqref{expex} with $\la\in(1,\infty)$. 
\eenu
\erem

The next theorem reports some basic properties of the solution found in Theorem~\ref{main3}:
\bthm\label{further} Let Assumption~\ref{Cass} be valid and $Y$ be the unique solution to Equation~\eqref{SPDE-var} in $B^{p,w}_{I,\loc}$. 
\benu
	\item For $(t,x),(\tau,\xi),(s,y)\in I\times\bbr^d$ define 	
	\begin{align} \tilde G(t,x;\tau,\xi;s,y)&:= \bigg(|G(t,x;s,y)-G(\tau,\xi;s,y)|^p\left(\int_\bbr |z|^p \,\pi(s,y,\dd z) + c(s,y)\right)\nonumber\\
	&\quad+ \big|[G(t,x;s,y)-G(\tau,\xi;s,y)]b_1(s,y)\big|\bone_{\{p\geq1\}}\bigg)w(s,y).\label{deftildeG2} \end{align}
If for all $(t,x)\in I\times\bbr^d$
	\beq\label{Lpcontcond} \int_I\int_{\bbr^d} \tilde G(t,x;\tau,\xi;s,y)\,\la(\dd s,\dd y) \to 0 \eeq
	whenever $(\tau,\xi)\to(t,x)$, then $Y$ is an $L^p$-continuous process, that is,
	\beq\label{LpY} \bbe[|Y(t,x)-Y(\tau,\xi)|^p] \to 0,\quad\text{whenever } (\tau,\xi)\to(t,x).\eeq
	\item Assume the case of Corollary~\ref{cor3} with $G(t,x;s,y)=g(t-s,x-y)$. Then \eqref{Lpcontcond} and therefore the conclusion of (1) hold automatically.
	\item $Y$ depends continuously on $Y_0$. In other words, if $Y$ and $Y^\prime$ are the solutions to \eqref{SPDE-var} with $Y_0, Y^\prime_0 \in B^{p,w}_{I,\loc}$ as force functions, respectively, then there exists a constant $C_{I,T,w}\in\bbr_+$ that may depend on $I$, $T$ and $w$, but is independent of $Y_0$, $Y^\prime_0$ such that
	\beq\label{contdep} \|Y-Y^\prime\|_{B^{p,w}_{I_T}} \leq C_{I,T,w} \|Y_0-Y_0^\prime\|_{B^{p,w}_{I_T}}. \eeq
	\eenu
\ethm

One of our basic motivations for studying Equation~\eqref{SPDE-var} on $I=\bbr$ is to construct stationary solutions. We show that if $G$ is of convolution form and $\La$ is homogeneous over space and time, then the stationarity of the solution in Theorem~\ref{main3} follows naturally.
\bthm\label{statsol}
Assume that $G(t,x;s,y)=g(t-s,x-y)$ and that $\La$ is a homogeneous L\'evy basis, satisfying the assumptions of Corollary~\ref{cor3}. Furthermore, suppose that for all $(t,x)\in\bbr_+\times\bbr^d$ we have
\beq\label{cag0} s\downarrow t,\quad y\downarrow x ~(\text{i.e. } y_i \downarrow x_i \text{ for all } i=1,\ldots,d) \quad\Longrightarrow\quad g(s,y)\to g(t,x),  \eeq
or that for all $(t,x)\in\bbr_+\times\bbr^d$ implication \eqref{cag0} holds with $\downarrow$ replaced by $\uparrow$.
  If $Y_0$ is $L^p$-continuous and jointly stationary with $\La$, then also $Y$ and $\La$ are jointly stationary.
\ethm

\bex[Continuation of Examples~\ref{ex1} and \ref{ex2}]\label{ex3} While the number $a$ in \eqref{heatkernel} did not play any role in Examples~\ref{ex1} and \ref{ex2}, this changes when we consider the stochastic heat equation on $I=\bbr$. Let $p\in(0,1+2/d)$ and set $p(d):=(1-p)d/2$. Then we have the following trichotomy: for $a>0$ we have
\beq\label{a1} \int_0^\infty \int_{\bbr^d} g_a^p(t,x)\,\dd(t,x) = (4\pi)^{p(d)}p^{-\frac{d}{2}}(ap)^{-1-p(d)}\Gamma(1+p(d)); \eeq
for $a=0$ we have for $T\in\bbr_+$
\beq\label{a0} \int_0^T \int_{\bbr^d} g_0^p(t,x)\,\dd(t,x) = \frac{{(4\pi)^{p(d)}}p^{-\frac{d}{2}}}{1+p(d)} T^{1+p(d)}, \eeq
which is of polynomial growth when $T\to\infty$; finally, if $a<0$, we have
\beq\label{1a} \int_0^T \int_{\bbr^d} g_a^p(t,x)\,\dd(t,x) = (4\pi)^{p(d)}p^{-\frac{d}{2}} \int_0^T \ee^{-apt} t^{p(d)}\,\dd t, \eeq
which grows faster than $\ee^{-apT}$ as $T\to\infty$. Thus, in the latter two cases, for Theorem~\ref{main3} to be applicable, the characteristics of $\La$ must decay fast enough at $-\infty$ to ensure the integrability conditions (4), (5) and (6) of Assumption~\ref{Cass}.

We will only focus on the case $a>0$. Given sufficiently strong decay properties of $\La$ at $-\infty$, the subsequent arguments can easily be transferred to the other two cases. First, we assume that $w\equiv1$ and that (1) and (2) of Assumption~\ref{Cass} hold. We further suppose the quasi-stationary case of \eqref{quasistat2}, and that the following conditions hold:
\beq\label{condcorr3} p<1+\frac{2}{d},\quad b_0\equiv0 \text{ if } p<1,\quad c\equiv0 \text{ if } p<2,\quad \zeta_p:=\int_\bbr |z|^p\,\pi_0(\dd z)<\infty. \eeq 
The only condition left is the size condition \eqref{intconst1} for $p\in(0,1)$ and \eqref{intconst2} for $p\in[1,2]$, respectively, before we can apply Corollary~\ref{cor3}. By \eqref{a1}, they are equivalent to
\beq\label{condex3} \zeta_p C_{\si,1}^p (4\pi)^{p(d)}p^{-\frac{d}{2}}(ap)^{-1-p(d)}\Gamma(1+p(d)) < 1 \eeq
in the case $p\in(0,1)$, and to
\beq\label{condex3-2} C_{\si,1}\left[C^\mathrm{BDG}_p \left((\zeta_p+\|c\|_{L^\infty_\bbr})(4\pi)^{p(d)}p^{-\frac{d}{2}}(ap)^{-1-p(d)}\Gamma(1+p(d))\right)^{1/p}+\|b_1\|_{L^\infty_\bbr}a^{-1} \right] < 1\eeq
in the case $p\in[1,2]$.

Finally, we would like to demonstrate how weighted spaces can be useful in Theorem~\ref{main3}. Let $a>0$ and $p\in(0,1+2/d)$ as before and define $w(t,x):=\ee^{\eta t}$ with $\eta\in\bbr$ satisfying $ap+\eta>0$. Assume that $Y_0\in B^{p,w}_{\bbr,\loc}$ and, if $\eta<0$, that $\si(0)=0$. Since
\begin{align*} &~\sup_{(t,x)\in\bbr\times\bbr^d} \int_{-\infty}^t\int_{\bbr^d} w^{-1}(t,x)g_a^p(t-s,x-y)w(s,y) \,\dd(s,y) = \int_0^\infty\int_{\bbr^d} g^p_a(s,y)\ee^{-\eta s}\,\dd(s,y)\\
=&~(4\pi)^{p(d)}p^{-\frac{d}{2}}(ap+\eta)^{-1-p(d)}\Gamma(1+p(d)),
\end{align*}
we have that in the conditions \eqref{condex3} and \eqref{condex3-2}, the term $ap$ is now replaced by $ap+\eta$. We draw two conclusions: if $Y_0$ is sufficiently small at $-\infty$, meaning $Y_0\in B^{p,w}_{\bbr,\loc}$ for some $\eta>0$, then the conditions \eqref{condex3} and \eqref{condex3-2} can be relaxed by using $ap+\eta$ instead of $ap$. Contrarily, if $\si(0)=0$, $\eta<0$, and the left-hand side of \eqref{condex3} or \eqref{condex3-2}, respectively, remains smaller than $1$ with $ap+\eta$ instead of $ap$, then one can even construct solutions with $Y_0\in B^{p,w}_{\bbr,\loc}$ that diverges at $-\infty$.
\eex

\section{Asymptotic stability}\label{assst}

In Theorems~\ref{main1}, \ref{main2} and \ref{main3} we have established solutions to \eqref{SPDE-var} that belong to the space $B^{p,w}_{I,\loc}$. In this section we will give criteria under which they even belong to the space $B^{p,w}_I$. Our primary focus is on the case where $\sup I=+\infty$, that is, we want to investigate whether solutions to \eqref{SPDE-var} are asymptotically $L^p$-stable. Moreover, we shall replace the Lipschitz condition on $\si$, which was essential in Sections~\ref{localex} and \ref{globalex}, by another growth condition, which, as we shall see, will determine the asymptotic behaviour of the solution. Of course, due to the possible non-Lipschitzianity of $\si$, we now have to \emph{assume} the existence of a solution in $B^{p,w}_{I,\loc}$. In fact, this approach  allows us to include solutions to \eqref{SPDE-var} with non-Lipschitz $\si$ which go beyond the results of the Sections~\ref{localex} and \ref{globalex} but are, for instance, studied in \citep{Mueller98, Mytnik02}.

Let us again start with a deterministic example that highlights the main features of the behaviour at infinity.
\bex\label{inftyex} Let $g\in L^1_{[0,\infty)}$, $f\in L^\infty_{[0,\infty)}$ and $v\in L^\infty_{[0,\infty),\loc}$ be positive functions satisfying
\beq\label{galeq1} v(t) = f(t) + \int_0^t g(t-s)v^\ga(s)\,\dd s,\quad t\in\bbr_+, \eeq
with $\ga\in(0,1]$. The question is under what conditions we have $v\in L^\infty_{[0,\infty)}$. It turns out that there is a fundamental difference between the cases $\ga\in(0,1)$ and $\ga=1$. In the former case, we always have $v\in L^\infty_{[0,\infty)}$. In fact, if we denote the convolution on the right-hand side of \eqref{galeq1} by $(g\ast v^\ga)(t)$, iteration of \eqref{galeq1} yields
\[ v = f + g\ast v^\ga = f + g \ast (f + g \ast v^\ga)^\ga = f + g \ast (f+ g\ast (f+g\ast v^\ga)^\ga)^\ga = \ldots \]
Using Young's inequality, we obtain
\begin{align*} \|v\|_{L^\infty_{[0,T]}} &\leq \|f\|_{L^\infty_{[0,\infty)}} + \|g\|_{L^1_{[0,\infty)}} \|v\|_{L^\infty_{[0,T]}}^\ga \leq \|f\|_{L^\infty_{[0,\infty)}} + \|g\|_{L^1_{[0,\infty)}} (\|f\|_{L^\infty_{[0,\infty)}} + \|g\|_{L^1_{[0,\infty)}}\|v\|_{L^\infty_{[0,T]}}^\ga)^\ga\\
&\leq \|f\|_{L^\infty_{[0,\infty)}} + \|g\|_{L^1_{[0,\infty)}} (\|f\|_{L^\infty_{[0,\infty)}} + \|g\|_{L^1_{[0,\infty)}}(\|f\|_{L^\infty_{[0,\infty)}} + \|g\|_{L^1_{[0,\infty)}} \|v\|_{L^\infty_{[0,T]}}^\ga)^\ga)^\ga \leq \ldots,  \end{align*}
or, equivalently, for every $T\in[0,\infty)$ and $n\in\bbn$
\[ \|v\|_{L^\infty_{[0,T]}} \leq a_n(T),\quad\text{where}\quad a_1(T):=\|v\|_{L^\infty_{[0,T]}},\quad a_{n+1}(T):=\|f\|_{L^\infty_{[0,\infty)}} + \|g\|_{L^1_{[0,\infty)}} (a_n(T))^\ga. \]
By induction it can be shown that $0\leq a_n(T)\leq a\vee a_1(T)$, where $a$ is the unique solution in $(0,\infty)$ of the equation
\[a-\|f\|_{L^\infty_{[0,\infty)}}-\|g\|_{L^1_{[0,\infty)}} a^\ga=0.\]
Note that $a$ does not depend on $T$, so we conclude that $\limsup_{n\to\infty} a_n(T)\leq a$ and $\|v\|_{L^\infty_{[0,T]}} \leq a$ for all $T\in[0,\infty)$. Hence we have $v \in L^\infty_{[0,\infty)}$ with $\|v\|_{L^\infty_{[0,\infty)}} \leq a$.

The situation is totally different for $\ga=1$. Then \eqref{galeq1} becomes
\beq\label{ga1} v(t)=f(t) + \int_0^t g(t-s) v(s)\,\dd s,\quad t\in\bbr_+, \eeq
which is the well known \emph{renewal equation}. If $f \in L^\infty_{[0,\infty)}$, one can show under some technical assumptions that the unique solution $v$ to \eqref{ga1} exhibits the following behaviour: if $\|g\|_{L^1_{[0,\infty)}}<1$, we have $v\in L^\infty_{[0,\infty)}$; if $\|g\|_{L^1_{[0,\infty)}} = 1$, the boundedness of $v$ depends on whether $f\in L^1_{[0,\infty)}$ or not; if $\|g\|_{L^1_{[0,\infty)}}>1$, then $v(t)\to\infty$ exponentially fast as $t\to\infty$. For precise statements with the required assumptions, we refer to \citep[Chap.~V]{Asmussen03}, especially to the Theorems~V.4.3 and V.7.1 and Proposition~V.7.4. 

In summary, whereas locally bounded solutions to \eqref{galeq1} with $\ga\in(0,1)$ are automatically globally bounded as soon as $f\in L^\infty_{[0,\infty)}$ and $g\in L^1_{[0,\infty)}$, the behaviour of the solution to \eqref{ga1} at infinity strongly depends on the \emph{size} of $\|g\|_{L^1_{[0,\infty)}}$. For a formalization of this example see also Lemma~\ref{assympholder} for $\ga\in(0,1)$ and Lemma~\ref{comparison} for $\ga=1$.
\eex

For Equation~\eqref{SPDE-var} the precise requirements are the following:
\bass\label{Dass} Let $p\in(0,2]$, $I$ be an interval and $w\colon I\times\bbr^d\to\bbr$ be a weight function satisfying $\sup_{(t,x)\in I\times\bbr^d} w^{-1}(t,x) < \infty$. We assume the following hypotheses:
\benu
	\item $Y_0\in B^{p,w}_I$.
	\item $\si\colon \bbr\to\bbr$ satisfies $|\si(x)|\leq |\si(0)|+C_{\si,2}|x|^\ga$ for all $x\in\bbr$ with some $\ga\in(0,1]$.
	\item Either $c(t,x)=0$ for all $(t,x)\in I\times\bbr^d$, or we have $2\ga\leq p$ and
	\beq\label{D-C} \sup_{(t,x)\in I\times\bbr^d} \int_I \int_{\bbr^d} w^{-1}(t,x)|G(t,x;s,y)|^2 w(s,y) c(s,y)\,\la(\dd s,\dd y) < \infty. \eeq 
	\item There exists $q\in(0,2]$ with $p\leq q$ and $q\ga\leq p$ such that
	\beq\label{D-nu} \sup_{(t,x)\in I\times\bbr^d} \int_I \int_{\bbr^d} \int_\bbr w^{-1}(t,x) |G(t,x;s,y)z|^p_q w(s,y)\,\nu(\dd s,\dd y,\dd z) < \infty.\eeq
	\item If $p\geq1$, then $\nu$ satisfies \eqref{helppi1} and 
	\beq\label{D-B1} \sup_{(t,x)\in I\times\bbr^d} \int_I \int_{\bbr^d} w^{-1}(t,x) |G(t,x;s,y)b_1(s,y)| w(s,y)\,\la(\dd s,\dd y)<\infty, \eeq
	and \eqref{D-B1} also holds with $w\equiv1$; if $p<1$, then there exist $\al\in(-\infty,2], \beta\in[0,\infty)$ satisfying \eqref{driftconv}, \eqref{meanconv} (with $\bbr_+$ replaced by $I$) and $(\al\vee\beta)\ga\leq p$ such that
	\beq\label{D-B0} \sup_{(t,x)\in I\times\bbr^d} \int_I \int_{\bbr^d} (F_0(s,y)\vee F_1(s,y))|G(t,x;s,y)|^\al_\beta \,\la(\dd s,\dd y) < \infty. \eeq
	\item At least one of the following three cases occurs:
	\item[(6a)] We have $\ga<1$, $q\ga<p$, $2\ga<p$ if $c\not\equiv0$ and $(\al\vee\beta)\ga<p$ if $p<1$.
	\item[(6b)] We have $p\in[1,2]$, and if we define for $(t,x), (s,y)\in I\times\bbr^d$
	\begin{align} 
	\bar G^{D,1}(t,x;s,y) &:= 2^{p-1}\left(\int_I\int_{\bbr^d} |G(t,x;s,y)b_1(s,y)|\,\la(\dd s,\dd y)\right)^{p-1}|G(t,x;s,y)b_1(s,y)|,\nonumber\\
	\bar G^{D,2}(t,x;s,y) &:= 2(C^\mathrm{BDG}_p)^2|G(t,x;s,y)|^2c(s,y),\nonumber\\
	\bar G^{D,3}(t,x;s,y) &:= 2^{p-1}(C^\mathrm{BDG}_p)^p \int_\bbr |G(t,x;s,y)z|^p\bone_{\{|G(t,x;s,y)z|>1\}}\,\pi(s,y,\dd z),\nonumber \\
	\bar G^{D,4}(t,x;s,y) &:= 2^{q-1}(C^\mathrm{BDG}_p)^q \int_\bbr |G(t,x;s,y)z|^q \bone_{\{|G(t,x;s,y)z|\leq1\}}\,\pi(s,y,\dd z),\nonumber\\
	G^{D,l}(t,x;s,y)&:=w^{-1}(t,x)G^{D,l}(t,x;s,y)w(s,y),\quad l=1,2,3,4, \label{GD}
	\end{align}
	then there exists a partition of $I$ into pairwise disjoint intervals $I_1, \ldots, I_k$ such that
	\beq\label{D-6b} \sup_{(t,x)\in I\times\bbr^d} \sup_{j=1,\ldots,k} \sum_{l=1}^4 C_{\si,2}\left( \int_{I_j}\int_{\bbr^d} G^{D,l}(t,x;s,y)\,\la(\dd s,\dd y) \right)^{1/p} < 1. \eeq
	\item[(6c)] We have $p\in(0,1)$, and if we define for $(t,x), (s,y) \in I\times\bbr^d$
	\begin{align}
	G^{D,1}(t,x;s,y) &:= 2^{(\al\vee\beta\vee1)-1}(F_0(s,y)\vee F_1(s,y))|G(t,x;s,y)|^\al_\beta,\nonumber\\
	G^{D,2}(t,x;s,y) &:= 2^{p+1}|G(t,x;s,y)|^2c(s,y),\nonumber\\
	G^{D,3}(t,x;s,y) &:= 2^p 2^{(q\vee1)-1} \int_\bbr |G(t,x;s,y)z|^p_q\,\pi(s,y,\dd z), \label{GD-2}
	\end{align}
	and
	\beq\label{rl} r_1:=\al\vee\beta,\quad r_2:=2,\quad r_3:=1,\eeq
	then there exists a partition of $I$ into pairwise disjoint intervals $I_1, \ldots, I_k$ such that
	\beq\label{D-6c} \sup_{(t,x)\in I\times\bbr^d} \sup_{j=1,\ldots,k} \sum_{l=1}^3 C_{\si,2}^{r_l} \int_{I_j}\int_{\bbr^d} G^{D,l}(t,x;s,y)\,\la(\dd s,\dd y) < 1. \eeq
\eenu
\eass

\bthm\label{assstat}
Let Assumption~\ref{Dass} be valid. If Equation \eqref{SPDE-var} has a solution $Y\in B^{p,w}_{I,\loc}$, it automatically also belongs to $B^{p,w}_I$.
\ethm

For quasi-stationary $G$ and $\La$, there is no significant simplification of Assumption~\ref{Dass} possible. Thus, we directly move to an example study.
\bex[Continuation of Examples~\ref{ex1}, \ref{ex2} and \ref{ex3}]\label{ex4} Let $I=[0,\infty)$, $a=0$ and $w\equiv1$. We assume that $Y\in B^p_{[0,\infty),\loc}$ solves 
\[ Y(t,x)=Y_0(t,x)+\int_0^t\int_{\bbr^d} g_0(t-s,x-y)\si(Y(s,y))\,\La(\dd s,\dd y),\quad (t,x)\in \bbr_+\times\bbr^d, \]
where $Y_0$ is given by \eqref{initial} and $\si$ satisfies condition (2) of Assumption~\ref{Dass} with $\ga\in(0,1]$. We want to find conditions that guarantee $Y\in B^p_{[0,\infty)}$. Let us check the requirements of Assumption~\ref{Dass}. (1) and (2) are clear. For (3), (4) and (5), the key observation is the following: for $p,q\in(0,2]$
\beq\label{pqa0} \int_0^\infty\int_{\bbr^d} |g_0(s,y)|^p_q\,\dd(s,y) < \infty \iff p\in(0,1+2/d) \text{ and } q\in(1+2/d,2]. \eeq
As a consequence of the last condition, unless in trivial cases, the classical stochastic heat equation with $a=0$ will be asymptotically unstable in dimensions $1$ and $2$. Only in dimensions $d\geq3$ there is a chance for asymptotic stability. We pose the following conditions:
\begin{align} &\la(\dd t,\dd x) = \dd(t,x), \quad \pi(t,x,\dd z)\leq \pi_0(\dd z), \quad p\in(0,1+2/d), \quad q\in(1+2/d,2],\nonumber\\
 &q\ga\leq p, \quad c\equiv0, \quad b_1\equiv0 \text{ if } p\geq1,\quad \La \text{ is symmetric if } p<1, \quad \int_\bbr |z|^q_p\,\pi_0(\dd z) <\infty.   \label{condcorr4} \end{align}
We notice that $\ga=1$ is not possible, and that $\int_\bbr |z|^p\,\pi_0(\dd z)<\infty$ is no longer sufficient, but $\pi_0$ must have a moment structure that is strictly better than its variation structure. Moreover, $c$ must be $0$; if $p\geq1$, only $\La\in\calm$ is possible; and if $p<1$, $\La$ is required to have no drift and a symmetric L\'evy measure. All this is because $g_0$ is not $L^p$-integrable on $\bbr_+\times\bbr^d$ for any $p\in(0,2]$. One readily sees that \eqref{condcorr4} implies conditions (3), (4) and (5). So if (6a) holds, we obtain $Y\in B^p_{[0,\infty)}$. In the case of (6b) or (6c), again a size condition has to be verified, which is analogous to the calculations in Example~\ref{ex3}. We leave the details to the reader. Note that in this example we have $\ga<1$, and therefore (6b) or (6c) is only needed in rare situations. Finally, for $a>0$ we refer the reader to the calculations in Example~\ref{ex3} again which can be re-used. In particular, one can find conditions for asymptotic stability in dimensions $1$ and $2$ this time.
\eex

\section{A series of lemmata}\label{proofs}

This section contains several lemmata that will play a crucial role in proving the main theorems in Section~\ref{proofmain}. First, we investigate the stochastic integral mapping in Equation~\eqref{SPDE-var}: fix some $\phi_0\in\tilde\calp$ and define for a predictable process $\phi\in\tilde\calp$ the process $J(\phi)$ by
\beq\label{eqQ} J(\phi)(t,x) := \phi_0(t,x) + \int_I \int_{\bbr^d} G(t,x;s,y)\si(\phi(s,y))\,\La(\dd s,\dd y) \eeq
for all $(t,x)\in I\in\bbr^d$ for which the stochastic integral exists, and set $J(\phi)(t,x):=+\infty$ otherwise. The next lemma, which is of crucial importance for all main results in this paper, relates the moment structure of $J(\phi)$ to that of $\phi$. 

\blem\label{momin} Let $w\colon I\times\bbr^d\to\bbr$ be a weight function.
\benu
	\item Suppose that Assumption~\ref{Cass} holds with $p\in(0,2]$ and recall the definition of $G^{C,1}$ and $G^{C,2}$ in \eqref{deftildeG}. Then for all $\phi\in\tilde\calp$ and $(t,x)\in I\times\bbr^d$, we have
	\begin{align} &~\frac{\|J(\phi)(t,x)\|_{L^p}}{(w(t,x))^{1/(p\vee1)}} \leq \frac{\|\phi_0(t,x)\|_{L^p}}{(w(t,x))^{1/(p\vee1)}}\nonumber\\
	&~\quad + \sum_{l=1}^2 \left(\int_I \int_{\bbr^d}  \frac{G^{C,l}(t,x;s,y)}{C^p_{\si,1}}\left(\frac{|\si(0)|^{p\wedge1}+C_{\si,1}^{p\wedge1}\|\phi(s,y)\|_{L^p}}{(w(s,y))^{1/(p\vee1)}}\right)^{p\vee1}\,\la(\dd s,\dd y)\right)^{1/(p\vee1)},\label{momin1}\end{align}
	where in the case $C_{\si,1}=0$ we use the convention $0/0:=1$.
	\item Furthermore, still under Assumption~\ref{Cass}, we have for any $\phi_1, \phi_2\in\tilde\calp$ for which the right-hand side of \eqref{momin1} is finite that
	\begin{align} &~\frac{\|J(\phi_1)(t,x)-J(\phi_2)(t,x)\|_{L^p}}{(w(t,x))^{1/(p\vee1)}}\nonumber\\
	\leq&~ \sum_{l=1}^2 \left(\int_I \int_{\bbr^d}  G^{C,l}(t,x;s,y)\left(\frac{\|\phi_1(s,y)-\phi_2(s,y)\|_{L^p}}{(w(s,y))^{1/(p\vee1)}}\right)^{p\vee1}\,\la(\dd s,\dd y)\right)^{1/(p\vee1)}. \label{momin2} \end{align}
	\item Let Assumption~\ref{Bass} or Assumption~\ref{Dass} be valid with $p\in[1,2]$. In the first case let $I=[0,\infty)$ and $w\equiv1$. Then the following holds for all $\phi\in\tilde\calp$ and $(t,x)\in I\times\bbr^d$:
	\begin{align} \frac{\|J(\phi)(t,x)\|_{L^p}}{(w(t,x))^{1/p}} &\leq \frac{\|\phi_0(t,x)\|_{L^p}}{(w(t,x))^{1/p}} + \frac{2[1+|\si(0)|+C_{\si,2}]}{(w(t,x))^{1/p}} \nonumber \\
	&\quad + (|\si(0)|+C_{\si,2})\sum_{l=1}^4 \left(\int_I\int_{\bbr^d} G^{D,l}(t,x;s,y)(w(s,y))^{-1}\,\la(\dd s,\dd y\right)^{1/p}\nonumber\\
	&\quad + \sum_{l=1}^4 C_{\si,2} \left(\int_I\int_{\bbr^d} \frac{G^{D,l}(t,x;s,y)}{(w(s,y))^{1-\rho}}\left(\frac{\|\phi(s,y)\|_{L^p}}{(w(s,y))^{1/p}}\right)^{p\rho}\,\la(\dd s,\dd y) \right)^{1/p},\label{momin3}\end{align} 
	where $G^{D,l}$ is defined by \eqref{GD}, and $\rho$ can be chosen as $\rho=(q\vee 2 \bone_{\{c\not\equiv0\}})\ga/p$ or $\rho=1$.
	\item Let Assumption~\ref{Bass} or Assumption~\ref{Dass} be valid with $p\in(0,1)$. In the first case let $I=[0,\infty)$ and $w\equiv1$. Then for all $\phi\in\tilde\calp$ and $(t,x)\in I\times\bbr^d$
	\begin{align}
	 \frac{\|J(\phi)(t,x)\|_{L^p}}{w(t,x)} &\leq \frac{\|\phi_0(t,x)\|_{L^p}}{w(t,x)} + \frac{2^{p+1}+1}{w(t,x)} \nonumber \\
	 &\quad+ \sum_{l=1}^3 (|\si(0)|^{r_l}_0+C_{\si,2}^{r_l})\int_I\int_{\bbr^d} G^{D,l}(t,x;s,y)(w(s,y))^{-1}\,\la(\dd s,\dd y) \nonumber \\
	 &\quad+ \sum_{l=1}^3 C_{\si,2}^{r_l} \int_I\int_{\bbr^d} \frac{G^{D,l}(t,x;s,y)}{(w(s,y))^{1-\rho}} \left(\frac{\|\phi(s,y)\|_{L^p}}{w(s,y)}\right)^\rho\,\la(\dd s,\dd y).
	\label{momin4}
	\end{align}
	where $G^{D,l}$ and $r_l$ are given by \eqref{GD-2} and \eqref{rl}, and $\rho=(q\vee 2\bone_{\{c\not\equiv0\}} \vee \al\vee\beta)\ga/p$ or $\rho=1$.
\eenu
\elem

\bpr It suffices to prove the lemma for $w\equiv1$: the general case follows if we divide the equations \eqref{momin1}, \eqref{momin2} and \eqref{momin3} by $w^{1/(p\vee1)}$. Throughout the proof, $(t,x)\in I\times\bbr^d$ is fixed, and the abbreviations $\Phi(s,y):=G(t,x;s,y)[\si(\phi_1(s,y))-\si(\phi_2(s,y))]$ and $\Psi(s,y):=G(t,x;s,y)\si(\phi(s,y))$ are used. Moreover, in the numerous integrals below, we will often drop the integration variables and use the shorthand notations $\iint_t := \int_{I_t}\int_{\bbr^d}$ and $\iiint_t := \int_{I_t} \int_{\bbr^d} \int_\bbr$. 

a) We first prove (2) when $p\geq1$.  To this end, we decompose 
\begin{align} \La(\dd t,\dd x) &= \left[\La^\cc(\dd t,\dd x) + \int_{\bbr} z\,(\mu-\nu)(\dd t,\dd x,\dd z) \right] +\left[B(\dd t,\dd x) + \int_{\bbr} z\bone_{\{|z|>1\}}\,\nu(\dd t,\dd x,\dd z)\right] \nonumber\\
&=: M(\dd t,\dd x)+B_1(\dd t,\dd x), \label{LamartB1}\end{align}
and obtain that $\|J(\phi_1)(t,x)-J(\phi_2)(t,x)\|_{L^p}$ is bounded by
\[ \|J^{(1)}(\phi_1)(t,x)-J^{(1)}(\phi_2)(t,x)\|_{L^p} + \|J^{(2)}(\phi_1)(t,x)-J^{(2)}(\phi_2)(t,x)\|_{L^p},\] 
where $J^{(1)}$ and $J^{(2)}$ are defined as in \eqref{eqQ} with $\La$ replaced by $M$ and $B_1$, respectively. 
For the $J^{(2)}$-part, H\"older's inequality yields
\begin{align}
\|J^{(2)}(\phi_1)(t,x)-J^{(2)}(\phi_2)(t,x)\|_{L^p} &\leq C_{\si,1}\left[\left(\iint_t |G|\,\dd |B_1|\right)^{p-1} \iint_t |G|\bbe[|\phi_1-\phi_2|^p] \,\dd |B_1|\right]^{1/p}\nonumber\\
&=\left(\iint_t G^{C,2}(t,x;s,y)\|\phi_1(s,y)-\phi_2(s,y)\|_{L^p}^p\,\la(\dd s,\dd y)\right)^{1/p}. \label{Bpart} \end{align}
For the $J^{(1)}$-part, we assume for the moment that the process
\beq\label{Nmart} N_\tau:=\iint_\tau G(t,x;s,y)[\si(\phi_1(s,y))-\si(\phi_2(s,y))]\,M(\dd s,\dd y)=\Phi\cdot M_\tau,\quad \tau\in I, \eeq
which is well defined by assumption, is a local martingale. Then we have by Definition~\ref{BDGineq} and the assumption that $c\equiv0$ for $p<2$
\begin{align} &~\|J^{(1)}(\phi_1)(t,x)-J^{(1)}(\phi_2)(t,x)\|_{L^p} \nonumber\\
\leq&~ C^{\mathrm{BDG}}_p \big\|[N]^{1/2}_t\big\|_{L^p} = C^{\mathrm{BDG}}_p \left\|\left(\iiint_t |\Phi z|^2\,\dd\mu + \iint_t |\Phi|^2 \,\dd C\right)^{1/2}\right\|_{L^p} \nonumber\\
\leq&~C^{\mathrm{BDG}}_p \bbe\left[\iiint_t |\Phi z|^p\,\dd\mu + \iint_t |\Phi|^2\,\dd C\right]^{1/p} = C^{\mathrm{BDG}}_p \bbe\left[\iiint_t |\Phi z|^p\,\dd\nu + \iint_t |\Phi|^2 \,\dd C\right]^{1/p} \nonumber\\
\leq&~\left(\iint_t G^{C,1}(t,x;s,y) \|\phi_1(s,y)-\phi_2(s,y)\|_{L^p}^p \,\la(\dd s, \dd y)\right)^{1/p}. \label{martcase}
\end{align}
Equations \eqref{Bpart} and \eqref{martcase} together imply \eqref{momin2} for $p\in[1,2]$. It remains to discuss whether $N$ in \eqref{Nmart} is a local martingale. Without loss of generality, we may assume that the right-hand side of \eqref{martcase} is finite; otherwise \eqref{momin2} becomes trivial. Let $\eps>0$ and $H\in\tilde\calp$ be a bounded function satisfying $|H(\om,s,y)|\leq \eps|\Phi(\om,s,y)|$ pointwise for all $(\om,s,y)\in\Om\times I\times\bbr^d$. Then $H\cdot M$ is a martingale such that we have by the Burkholder-Davis-Gundy inequality
\begin{align} \sup_{\tau\in I}\|H\cdot M_\tau\|_{L^p} &\leq C^{\mathrm{BDG}}_p \left\|\left(\iiint_t |H z|^2\,\dd\mu + \iint_t |H|^2 \,\dd C\right)^{1/2}\right\|_{L^p} \nonumber\\
&\leq \eps C^{\mathrm{BDG}}_p \left\|\left(\iiint_t |\Phi z|^2\,\dd\mu + \iint_t |\Phi|^2 \,\dd C\right)^{1/2}\right\|_{L^p}. \label{semivar}\end{align}
The right-hand side of \eqref{semivar} is finite by \eqref{martcase}. Moreover, as $\eps\downarrow0$, it goes to $0$ independently of $H$. Thus, \citep[Prop.~4.9b]{Bichteler83} is applicable (the extension of this proposition to intervals $I$ different from $I=\bbr_+$ is straightforward) and shows that $N$ is indeed a local martingale. 

b) We prove (2) when $p<1$. By hypothesis, $\La$ is L\'evy basis without drift. Thus, using the basic estimate $\left|\sum_{i=1}^\infty a_i\right|^p \leq \sum_{i=1}^\infty |a_i|^p$, we obtain
\begin{align*} &~\|J(\phi_1)(t,x)-J(\phi_2)(t,x)\|_{L^p} = \left\|\iiint_t \Phi z\,\dd\mu\right\|_{L^p} \leq \bbe\left[\iiint_t |\Phi z|^p\,\dd\mu\right]\\
 \leq&~ C^p_{\si,1} \iiint_t |Gz|^p \bbe[|\phi_1-\phi_2|^p]\,\dd\nu = \iint_t G^{C,1}(t,x;s,y) \|\phi_1(s,y)-\phi_2(s,y)\|_{L^p}\,\la(\dd s,\dd y), \end{align*}
which is \eqref{momin2}.

c) Because the Lipschitz condition on $\si$ implies $|\si(x)|\leq |\si(0)|+C_{\si,1}|x|$ for all $x\in\bbr$, (1) can be deduced in complete analogy to a) and b).

d) We prove (3). To this end, we again consider the decomposition $\La=M+B_1$ as in \eqref{LamartB1}. Using Definition~\ref{BDGineq}, Jensen's inequality and the hypothesis that $q\ga\leq p$ and $2\ga\bone_{\{c\not\equiv0\}}\leq p$, we obtain
\allowdisplaybreaks
\begin{align} &~\|\Psi\cdot M_t\|_{L^p} \leq C^\mathrm{BDG}_p \left\|\left(\iiint_t |\Psi z|^2\,\dd\mu + \iint_t |\Psi|^2 \,\dd C\right)^{1/2}\right\|_{L^p} \nonumber\\
\leq&~ C^\mathrm{BDG}_p \left(\bbe\left[ \iiint_t |\Psi z|^q \bone_{\{|Gz|\leq1\}}\,\dd\nu \right]^{p/q} + \bbe\left[\iiint_t |\Psi z|^p \bone_{\{|Gz|>1\}}\,\dd\nu + \left(\iint_t |\Psi|^2 \,\dd C\right)^{p/2}\right]\right)^{1/p}\nonumber\\
\leq&~C^\mathrm{BDG}_p \Bigg[\left(2^{q-1}\iiint_t |Gz|^q \bone_{\{|Gz|\leq1\}} (|\si(0)|^q+C^q_{\si,2}\|\phi\|_{L^p}^{q\ga})\,\dd\nu \right)^{1/q}\nonumber \\
&\quad+\left(2^{p-1}\iiint_t |Gz|^p \bone_{\{|Gz|>1\}} (|\si(0)|^p+C^p_{\si,2}\|\phi\|_{L^p}^{p\ga})\,\dd\nu\right)^{1/p}\nonumber\\
&\quad+ \left(2\iint_t |G|^2 (|\si(0)|^2+C^2_{\si,2}\|\phi\|_{L^p}^{2\ga})\,\dd C\right)^{1/2}\Bigg]\nonumber\\
\leq&~C^\mathrm{BDG}_p \Bigg[\left(2^{q-1}\iiint_t |Gz|^q \bone_{\{|Gz|\leq1\}} (|\si(0)|^q+C^q_{\si,2}+C_{\si,2}^q\|\phi\|_{L^p}^{p\rho})\,\dd\nu \right)^{1/q}\nonumber \\
&\quad+\left(2^{p-1}\iiint_t |Gz|^p \bone_{\{|Gz|>1\}} (|\si(0)|^p+C^p_{\si,2}+C^p_{\si,2}\|\phi\|_{L^p}^{p\rho})\,\dd\nu\right)^{1/p}\nonumber\\
&\quad+ \left(2\iint_t |G|^2 (|\si(0)|^2+C^2_{\si,2}+C^2_{\si,2}\|\phi\|_{L^p}^{p\rho})\,\dd C\right)^{1/2}\Bigg]\nonumber\\
\leq& (|\si(0)|+C_{\si,2}) \left[2+\sum_{l=2}^4 \left(\int_I\int_{\bbr^d} G^{D,l}(t,x;s,y)\,\la(\dd s,\dd y)\right)^{1/p}\right] \nonumber\\
&\quad+ 2 + \sum_{l=2}^4 C_{\si,2} \left(\int_I\int_{\bbr^d} G^{D,l}(t,x;s,y)\|\phi(s,y)\|_{L^p}^{p\rho}\,\la(\dd s,\dd y)\right)^{1/p}. \label{mart3}
\end{align}
\allowdisplaybreaks[0]
Again, one can justify that $\Psi\cdot\La$ is indeed a well defined local martingale whenever the right-hand side of \eqref{momin3} is finite.
For the $B_1$-integral, another application of H\"older's inequality demonstrates
\begin{align} \|(\Psi\cdot B_1)_t\|_{L^p} &\leq \left[2^{p-1}\left(\iint_t |G|\,\dd|B_1|\right)^{p-1}\iint_t |G|(|\si(0)|^p + C_{\si,2}^p\|\phi\|_{L^p}^{p\ga})\,\dd|B_1|\right]^{1/p}\nonumber\\
&\leq \left[2^{p-1}\left(\iint_t |G|\,\dd|B_1|\right)^{p-1}\iint_t |G|(|\si(0)|^p + C_{\si,2}^p+C_{\si,2}^p\|\phi\|_{L^p}^{p\rho})\,\dd|B_1|\right]^{1/p}. \label{B3} \end{align}
Equation~\eqref{momin3} now follows from \eqref{mart3} and \eqref{B3}.

e) We consider the last part (4). In this case we directly use the canonical decomposition of $\Psi\cdot\La$:
\[ \Psi\cdot\La_t = \Psi\cdot\La^\cc_t + \iiint_t \Psi z\bone_{\{|\Psi z|\leq1\}}\,\dd(\mu-\nu) + \iiint_t \Psi z\bone_{\{|\Psi z|>1\}}\,\dd\mu + B^{\Psi\cdot\La}_t =: J^1 + J^2 + J^3 + J^4, \]
where
\[ B^{\Psi\cdot\La}(\dd t,\dd x) = \Psi(t,x)\left[b(t,x)+ \int_\bbr z(\bone_{\{|\Psi(t,x)z|\leq1\}}-\bone_{\{|z|\leq1\}})\,\pi(t, x,\dd z)\right]\,\la(\dd t,\dd x).\]
We begin with $J^1$:
\begin{align*} \|J^1\|_{L^p} &\leq (C^{\mathrm{BDG}}_1)^p \bbe\left[\left(\iint_t |\Psi|^2 \,\dd C\right)^{1/2}\right]^p \leq 2^p \left(\iint_t \bbe[|\Psi|^2] \,\dd C\right)^{p/2} \\
&\leq 2^p \left(1+2\iint_t G^2 (|\si(0)|^2+C_{\si,2}^2\|\phi\|_{L^p}^{2\ga/p})\,\dd C\right)\\
&\leq 2^p \left(1+2\iint_t G^2 (|\si(0)|^2+C_{\si,2}^2+C_{\si,2}^2\|\phi\|_{L^p}^\rho)\,\dd C\right).  \end{align*}
For the jumps part, we obtain
\begin{align*} &~\|J^2+J^3\|_{L^p} \leq (C^{\mathrm{BDG}}_1)^p \bbe\left[\left(\iiint_t |\Psi z|^2 \bone_{\{|\Psi z|\leq1\}}\,\dd\mu\right)^{1/2}\right]^p + \bbe\left[\iiint_t|\Psi z|^p\bone_{\{|\Psi z|>1\}}\,\dd\mu\right] \\
\leq&~ 2^p \left(\iiint_t \bbe[|\Psi z|^q\bone_{\{|\Psi z|\leq1\}}] \,\dd\nu\right)^{p/2} + \iiint_t \bbe[|\Psi z|^p\bone_{\{|\Psi z|>1\}}]\,\dd\nu \\
\leq&~ 2^p\left(1+\iiint_t \bbe[|\Psi z|^p_q]\,\dd\nu\right) \leq 2^p \left(1+2^{(q\vee1)-1}\iiint_t |Gz|^p_q (|\si(0)|^q_0+C_{\si,2}^q\|\phi\|_{L^p}^{q\ga/p})\,\dd\nu\right)\\
\leq&~ 2^p \left(1+2^{(q\vee1)-1}\iiint_t |Gz|^p_q (|\si(0)|^q_0+C_{\si,2}^q+C_{\si,2}^q\|\phi\|_{L^p}^\rho)\,\dd\nu\right).  \end{align*}
%
Finally, since
\[ |J^4| \leq \iint_t |\Psi(s,y)| \left|b(s,y) + \int_\bbr \left[z\bone_{\{|z|\in (1,|\Psi(s,y)|^{-1}]\}} - z\bone_{\{|z|\in(|\Psi(s,y)|^{-1},1]\}}\right]\,\pi(s, y,\dd z)\right|\,\la(\dd s,\dd y), \]
we deduce the following bound for $J^4$ from Assumption~\ref{Bass}(9) or Assumption~\ref{Dass}(4), respectively:
\begin{align*}
\|J^4\|_{L^p}&\leq \bbe\left[\iint_t |\Psi|(|\Psi|^{\beta-1} F_1 \bone_{\{|\Psi|\leq1\}} + |\Psi|^{\al-1} F_0\bone_{\{|\Psi|>1\}})\,\dd\la\right]^p \leq \left(\iint_t (F_0\vee F_1)\bbe[|\Psi|^\al_\beta]\,\dd\la\right)^p\\
&\leq 1+2^{(\al\vee\beta\vee1)-1} \iint_t (F_0\vee F_1) |G|^\al_\beta \left(|\si(0)|^{\al\vee\beta}_0+C_{\si,2}^{\al\vee\beta}\|\phi\|_{L^p}^{(\al\vee\beta)\ga/p}\right)\,\dd\la\\
&\leq 1+2^{(\al\vee\beta\vee1)-1} \iint_t (F_0\vee F_1) |G|^\al_\beta \left(|\si(0)|^{\al\vee\beta}_0+C_{\si,2}^{\al\vee\beta}+C_{\si,2}^{\al\vee\beta}\|\phi\|_{L^p}^\rho\right)\,\dd\la.
\end{align*}
In combination with the estimates for $J^1$, $J^2$ and $J^3$, this finishes the proof of \eqref{momin4}.
\epr

The next lemma allows us to take good versions of the stochastic integral process \eqref{eqQ}:
\blem\label{predver}
For every $\phi\in\tilde\calp$ there exists a predictable modification of $J(\phi)$, that is, a $(-\infty,\infty]$-valued process $\bar J(\phi)\in\tilde\calp$ such that for each $(t,x)\in I\times\bbr^d$ we have $J(\phi)(t,x)=\bar J(\phi)(t,x)$ a.s.
\elem
\bpr
The set $A$ of all $(t,x)\in I\times\bbr^d$ for which $G(t,x;\cdot,\cdot)\si(\phi)$ is integrable w.r.t. $\La$ is deterministic by definition, and by \citep[Thm.~4.1]{Chong14} and Fubini's theorem also measurable. It follows that there exists a measurable modification $J^\mathrm{m}(\phi)$ of $J(\phi)$: set $J^\mathrm{m}(\phi)=\infty$ on $A^\Comp$ and use \citep[Thm.~1]{Lebedev95} on $A$. Next, define $\ppp J(\phi)$ as the extended predictable projection of $J^\mathrm{m}(\phi)$ in the sense of \citep[Thm.~I.2.28]{Jacod03}. By \citep[Prop.~3]{Stricker78} we may choose $\ppp J(\phi)(t,x)$ measurably in $x$. And indeed, $\ppp J(\phi)$ is still a modification of $J(\phi)$ since for each $(t,x)\in I\times\bbr^d$ we have a.s.
\[ \ppp J(\phi)(t,x) = \bbe[J^\mathrm{m}(\phi)(t,x)\,|\,\calf_{t-}] = \int_{I_t} \int_{\bbr^d} G(t,x;s,y)\si(\phi(s,y))\,\La(\dd s, \dd y) = J(\phi)(t,x).\]
\epr

We proceed with a discretization result for stochastic integrals:
\blem\label{Riemann} Let $I\subseteq\bbr$ be an interval and $w\equiv1$, and assume that $G$, $\si$ and $\La$ satisfy (2)--(6) of Assumption~\ref{Cass}. Fix some $(t,x)\in I\times\bbr^d$ and assume that $G(t,x;\cdot,\cdot)$ has the following properties: for all $(s,y)\in I_t,\times\bbr^d$ we have
\beq\label{cag} r\uparrow s,\quad z\uparrow y ~(\text{i.e. } z_i \uparrow y_i \text{ for all } i=1,\ldots,d) \quad\Longrightarrow\quad G(t,x;r,z)\to G(t,x;s,y),  \eeq
and for some $\eps>0$ the function $G^\ast_\eps(t,x;s,y):= \sup_{r\in I, s-\eps<r\leq s, |y-z|<\eps} |G(t,x;r,z)|$ satisfies
\begin{align} &~\int_{I_t} \int_{\bbr^d} \Bigg(|G^\ast_\eps(t,x;s,y)|^p \left(\int_{\bbr} |z|^p\,\pi(s,y,\dd z) + c(s,y)\right) \nonumber\\ &~\quad+|G^\ast_\eps(t,x;s,y)b_1(s,y)|\bone_{\{p\geq1\}}\Bigg)\,\la(\dd s,\dd y) < \infty.\label{supbounded}\end{align}
Moreover, we specify discretization schemes for both time and space: first, we choose for each $N\in\bbn$ a number $k(N)\in\bbn\cup\{\infty\}$ of time points $(s_i^N)_{i=1}^{k(N)}\subseteq I_t$ such that
\[ s_i^N < s_{i+1}^N,\quad\text{and}\quad s_1^N \downarrow \inf I,\quad s_{k(N)}^N \uparrow t, \quad \sup_{i=1,\ldots, k(N)-1} |s_{i+1}^N-s_i^N| \downarrow 0 \quad\text{as } N\uparrow\infty; \]
and second, we fix for each $N\in\bbn$ a number $l(N)\in\bbn\cup\{\infty\}$ of non-empty pairwise disjoint hyperrectangles $(Q_j^N=(a_j^N,b_j^N])_{j=1}^{l(N)} \subseteq \bbr^d$ satisfying
\[ \bigcup_{j=1}^{l(N)} Q_j^N \uparrow \bbr^d\quad\text{and}\quad \sup_{j=1,\ldots,l(N)} \diam(Q_j^N) \downarrow 0\quad\text{as } N\uparrow\infty. \]
\benu
\item If $\phi\in B^p_{I,\loc}$ is an $L^p$-continuous process (cf. \eqref{LpY}), then the stochastic integral $J(\phi)(t,x)$ is well defined and
\beq\label{approx} \phi_0(t,x) + \sum_{i=1}^{k(N)-1}\sum_{j=1}^{l(N)} G(t;x;s_i^N,a_j^N)\si(\phi(s_i^N,a_j^N))\La\big((s_i^N,s_{i+1}^N] \times Q_j^N \big) \to J(\phi)(t,x) \eeq
in $L^p$ as $N\to\infty$.
\item The statement of (1) remains true if we replace $\uparrow$ in \eqref{cag} by $\downarrow$, and at the same time replace $G(t,x;s_i^N,a_j^N)$ by $G(t,x;s_{i+1}^N,b_j^N)$ in \eqref{approx}. 
\eenu
\elem
\bpr Part (2) is proved in the same fashion as part (1). That the stochastic integral $J(\phi)(t,x)$ exists, is a consequence of Lemma~\ref{momin}(1), the assumptions posed on $G$ and $\La$, and the fact that $\phi\in B^p_{I,\loc}$. To prove \eqref{approx}, let us call its left-hand side $J^N(\phi)(t,x)$. It follows that
\begin{align*} J^N(\phi)(t,x)&=\phi_0(t,x)+\int_{I_t}\int_{\bbr^d} H^N(t,x;s,y)\,\La(\dd s,\dd y), \quad\text{where}\\
 H^N(t,x;s,y) &=  \sum_{i=1}^{k(N)-1}\sum_{j=1}^{l(N)} G(t;x;s_i^N,a_j^N)\si(\phi(s_i^N,a_j^N)) \bone_{(s_i^N,s_{i+1}^N] \times Q_j^N}(s,y). \end{align*}
We notice that $H^N(t,x;s,y)=0$ if $(s,y)$ does not belong to the set $A^N := (s_1^N,s_{k(N)}^N]\times \bigcup_{j=1}^{l(N)} Q_j^N$, and that for each $(s,y)\in(\inf I,t)\times\bbr^d$ we have $\bone_{(A^N)^\Comp}(s,y) \to 0$ as $N\to\infty$. Now, we distinguish between two cases: first, if $p<1$, or $p\geq1$ and $\La\in\calm$, then similar calculations as done for Lemma~\ref{momin}(2) lead to (set $H(t,x;s,y):=G(t,x;s,y)\si(\phi(s,y))$)
\begin{align} &~\bbe[|J(\phi)(t,x)-J^N(\phi)(t,x)|^p]\nonumber\\
\leq&~ (C^\mathrm{BDG}_p)^p \int_{I_t}\int_{\bbr^d} \bbe[|H(t,x;s,y)-H^N(t,x;s,y)|^p] \left(\int_{\bbr} |z|^p\,\pi(s,y,\dd z) + c(s,y)\right)\,\la(\dd s,\dd y) \nonumber\\
\leq&~ \int_{(A^N)^\Comp} \frac{G^{C,1}(t,x;s,y)}{C_{\si,1}^p} \bbe[|\si(\phi(s,y))|^p]\,\la(\dd s,\dd y) \nonumber\\
&~\quad+ \iint_{A^N} \frac{G^{C,1}(t,x;s,y)}{C_{\si,1}^p} \sum_{i,j} \bbe[|\si(\phi(s,y))-\si(\phi(s_i^N,a_j^N))|^p] \bone_{(s_i^N,s_{i+1}^N] \times Q_j^N}(s,y) \,\la(\dd s,\dd y) \nonumber\\
&~\quad+ (C^\mathrm{BDG}_p)^p\iint_{A^N} \sum_{i,j} |G(t,x;s,y)-G(t,x;s_i^N,a_j^N)|^p \bbe[|\si(\phi(s_i^N,a_j^N))|^p] \bone_{(s_i^N,s_{i+1}^N] \times Q_j^N}(s,y)\nonumber\\
&~\quad\cdot\left(\int_{\bbr} |z|^p\,\pi(s,y,\dd z) + c(s,y)\right)\,\la(\dd s,\dd y) \nonumber\\
=&\!\!: I^N_1 + I^N_2 + I^N_3. \label{I3}
\end{align}
Since $\phi\in B^p_{I,\loc}$ and $G^{C,1}$ is integrable w.r.t. $\la$ by hypothesis, $I^N_1\to0$ as $N\to\infty$ by dominated convergence. Next, as a consequence of the $L^p$-continuity of $\phi$ and the refining properties of our discretization scheme, the sum within $I^N_2$ goes to $0$ pointwise for each $(s,y)\in I_t\times\bbr^d$. Moreover, this sum is majorized by $2\|\si(\phi)\|_{B^p_{I_t}}$ such that also $I^N_2\to0$ as $N\to\infty$. Regarding $I^N_3$, we obtain as an upper bound
\begin{align*} I^N_3&\leq (C^\mathrm{BDG}_p)^p\|\si(\phi)\|_{B^p_{I_t}} \iint_{A^N} \left|G(t,x;s,y)-\sum_{i,j}G(t,x;s_i^N,a_j^N)\bone_{(s_i^N,s_{i+1}^N] \times Q_j^N}(s,y) \right|^p \\ 
&\quad\cdot\left(\int_{\bbr} |z|^p\,\pi(s,y,\dd z) + c(s,y)\bone_{\{p=2\}}\right)\,\la(\dd s,\dd y).
\end{align*}
Because of \eqref{cag}, the integrand in the last line goes to $0$ as $N\to\infty$, pointwise for $(s,y)\in I_t\times\bbr^d$. Moreover, it is dominated by $2G^\ast_\eps$, when $\eps$ is chosen according to \eqref{supbounded} and $N$ is large enough such that $\sup_{i=1,\ldots, k(N)-1} |s_{i+1}^N-s_i^N|$ and $\sup_{j=1,\ldots,l(N)} \diam(Q_j^N)$ are smaller than $\eps$. By dominated convergence, we conclude $I^N_3\to0$ as $N\to\infty$. 

It remains to discuss the case $p\geq1$ and $\La\notin\calm$. As in Lemma~\ref{momin}(2), we decompose $\La=M+B_1$, where $M$ is a martingale measure and $B_1$ the drift measure. For $M$ we can apply the calculations above. For $B_1$ we obtain an analogous decomposition as in \eqref{I3}: $G^{C,1}$ is replaced by $G^{C,2}$, and instead of the Burkholder-Davis-Gundy constants, the factor 
\[ \left(\int_{I_t}\int_{\bbr^d} \sum_{i,j} |G(t,x;s,y)-G(t,x;s_i^N,a_j^N)| \bone_{(s_i^N,s_{i+1}^N] \times Q_j^N}(s,y)\,|B_1|(\dd s,\dd y)\right)^{p-1} \]
appears. But this also goes to $0$ as $N\to\infty$, as desired.
\epr

The next lemma concerns the solvability of deterministic integral equations and provides a comparison result. Certainly, there is a huge literature on deterministic Volterra equations, but we did not find a reference completely satisfying our purposes. Thus, we decided to include the proof, which is also very instructive for the proofs of the main theorems below. 

\blem\label{comparison}
Let $I\subseteq\bbr$ be an interval and $\la$ a positive measure on $(I\times \bbr^d,\calb(I\times\bbr^d))$ and $p\in [1,\infty)$. Further suppose that for every $l\in\bbn$ we have a positive measurable function $G^{(l)}\colon (I\times \bbr^d)^2 \to \bbr$ with $G^{(l)}(t,\cdot;s,\cdot)\equiv 0$ for $s>t$. Moreover, assume that there exists $k\in\bbn$ and a partition of $I$ into pairwise disjoint intervals $I_1,\ldots, I_k$ such that 
\beq\label{partition} \rho:=\sup_{(t,x)\in I\times \bbr^d} \sup_{j=1,\ldots,k} \sum_{l=1}^\infty \left(\int_{I_j}\int_{\bbr^d} G^{(l)}(t,x;s,y)\,\la(\dd s,\dd y)\right)^{1/{p}}<1. \eeq

Then the following statements hold:
\benu
\item Let $(v_n)_{n\in\bbn}$ be a sequence of positive functions in $L^{\infty}_I$ satisfying
\beq\label{itv} v_{n+1}(t,x) \leq \sum_{l=1}^\infty \left(\int_I\int_{\bbr^d} G^{(l)}(t,x;s,y) (v_n(s,y))^{p}\,\la(\dd s,\dd y)\right)^{1/{p}},\quad n\in\bbn. \eeq
Then $\sum_{n=1}^\infty \|v_n\|_{L^{\infty}_I}$ is finite. In particular, $v_n$ converges in $L^{\infty}_I$ to $0$.
\item For every positive $f\in L^{\infty}_I$ the equation 
\beq\label{volterra} v(t,x) = f(t,x) + \sum_{l=1}^\infty \left(\int_I \int_{\bbr^d} G^{(l)}(t,x;s,y) (v(s,y))^{p} \,\la(\dd s,\dd y)\right)^{1/{p}},\quad (t,x)\in I\times \bbr^d, \eeq
has a unique solution $v\in L^{\infty}_I$. Furthermore, this solution $v$ is positive.
\item If $\bar v\in L^{\infty}_I$ is a positive function satisfying
\beq\label{compare} \bar v(t,x)\leq f(t,x) + \sum_{l=1}^\infty \left(\int_I \int_{\bbr^d} G(t,x;s,y) (\bar v(s,y))^{p}\,\la(\dd s,\dd y)\right)^{1/{p}},\quad (t,x)\in I\times \bbr^d, \eeq
then we have $\bar v(t,x) \leq v(t,x)$ for all $(t,x)\in I\times \bbr^d$. In particular, if $f\equiv0$, then $v\equiv\bar v\equiv0$.
\eenu
\elem

\bpr a) We start with (1). Let $I=I_1\cup\ldots\cup I_k$ be as in the hypothesis and suppose that the intervals $I_j$ are arranged in increasing order (i.e. $\sup I_j = \inf I_{j+1}$). Furthermore, define for $\phi\in L^\infty_I$ and $(t,x)\in I\times \bbr^d$
\begin{align} \|\phi\|_{G^{(l)},p}(t,x) &:= \left(\int_I\int_{\bbr^d} G^{(l)}(t,x;s,y)|\phi(s,y)|^{p}\,\la(\dd s,\dd y)\right)^{1/{p}} \nonumber \\
\|\phi\|_{G^{(l)},p,j}(t,x) &:= \left(\int_{I_j}\int_{\bbr^d} G^{(l)}(t,x;s,y)|\phi(s,y)|^{p}\,\la(\dd s,\dd y)\right)^{1/{p}},\quad l\in\bbn,\quad j=1,\ldots,k. \label{Gnorms}\end{align}
Obviously, we have $\|\phi\|_{G^{(l)},p}(t,x) \leq \sum_{j=1}^k \|\phi\|_{G^{(l)},p,j}(t,x)$ for each $(t,x)\in I\times \bbr^d$ and $l\in\bbn$. Hence, it follows from \eqref{itv} that
\beq\label{vit} v_{n+1} \leq \sum_{l=1}^\infty \|v_n\|_{G^{(l)},p} \leq \sum_{j=1}^k \sum_{l=1}^\infty \|v_n\|_{G^{(l)},p,j}, \eeq
an equation that holds pointwise for all $(t,x)\in I\times \bbr^d$ and for all $n\in\bbn$. Iterating \eqref{vit} $n$ times, together with the subadditivity of the functional $\|\cdot\|_{G^{(l)},p,j}$, yields
\begin{align} v_{n+1} &\leq \sum_{j_1=1}^k \sum_{l_1=1}^\infty \|v_n\|_{G^{(l_1)},p,j_1} \leq \sum_{j_1,j_2=1}^k \sum_{l_1=1}^\infty \left\| \sum_{l_2=1}^\infty \|v_{n-1}\|_{G^{(l_2)},p,j_2} \right\|_{G^{(l_1)},p,j_1} \leq \ldots \nonumber\\
&\leq \sum_{j_1,\ldots,j_n=1}^k  \sum_{l_1=1}^\infty \left\| \sum_{l_2=1}^\infty \left\| \ldots \sum_{l_n=1}^\infty \|v_1\|_{G^{(l_n)},p,j_n} \ldots \right\|_{G^{(l_2)},p,j_2} \right\|_{G^{(l_1)},p,j_1}\label{vpin0}. \end{align}
Observe that the Volterra property of $G$ implies that on the right-hand side of \eqref{vpin0}, only those summands are non-zero for which $j_1\geq\ldots\geq j_n$. Since there are exactly $\binom{n+k-1}{n}$ such sequences, and $\sup_{j=1,\ldots,k} \left\| \sum_{l=1}^\infty \|1\|_{G^{(l)},p,j}\right\|_{L^\infty_I} = \rho$, we deduce that
\beq\label{summable} \sum_{n=1}^\infty \|v_n\|_{L^{\infty}_I} \leq \|v_1\|_{L^{\infty}_I} \sum_{n=0}^\infty \binom{n+k-1}{n} \rho^n < \infty \eeq
by the ratio test and the fact that $\rho<1$. 

b) Next we prove (2) and construct a solution to \eqref{volterra} by Picard iteration. Define $v^0(t,x)=f(t,x)$ and for $n\in\bbn$,
\beq\label{picit} v^n(t,x):=f(t,x) + \sum_{l=1}^\infty \left(\int_I\int_{\bbr^d} G^{(l)}(t,x;s,y) (v^{n-1}(s,y))^{p}\,\la(\dd s,\dd y)\right)^{1/{p}},\quad (t,x)\in I\times \bbr^d. \eeq
Since $G$ satisfies \eqref{partition}, $f$ belongs to $L^{\infty}_I$, and both functions are positive, $v^n$ is by induction again a positive function in $L^{\infty}_I$. Now form the difference sequence $u^n:=|v^{n+1}-v^n|$ for $n\in\bbn$, which satisfies property \eqref{itv} by the reverse triangle inequality. By (1), $\sum_{n=1}^\infty \|u^n\|_{L^{\infty}_I} < \infty$, in other words, $v$ as the limit in $L^{\infty}_I$ of $v^n$ exists. Of course, $v$ is positive. Moreover, taking the limit on both sides of \eqref{picit}, we conclude that $v$ indeed satisfies \eqref{volterra}. The uniqueness part follows by applying part (1) to the difference of two solutions in $L^{\infty}_I$.

c) For $\phi \in L^\infty_I$ set $I_{f}(\phi):=f+\sum_{l=1}^\infty \|\phi\|_{G^{(l)},p}$, which again belongs to $L^\infty_I$. By \eqref{picit}, we have $v^n = I^{(n)}_{f}(f)$, which is the $n$-fold iteration $I_f(I_f(\ldots I_f(f)\ldots))$. Moreover, by \eqref{compare}, 
\[ \bar v \leq I_{f}(\bar v) \leq I_{f}\big(I_{f}(\bar v)\big) \leq \ldots \leq I^{(n)}_{f}(\bar v) \leq I^{(n-1)}_{f}(f) + I^{(n)}_{0}(\bar v) = v^{n-1} + I^{(n)}_{0}(\bar v). \]
As shown in a), $v^{n-1}$ converges to $v$ uniformly on $I\times \bbr^d$. In addition, $I^{(n)}_{0}(\bar v)$ is less or equal to the right-hand side of \eqref{vpin0} when $v_1$ is replaced by $\bar v$. Thus, the considerations in a) show that $I^{(n)}_{0}(\bar v) \leq \|\bar v\|_{L^\infty_I} \binom{n+k-1}{n} \rho^n \to 0$ as $n\to\infty$, which implies (3).
\epr

The next lemma concerns the asymptotic behaviour of deterministic Volterra equations with a fractional nonlinearity:
\blem\label{assympholder} Let $I$, $p$ and $G^{(l)}$ be as in Lemma~\ref{comparison}. Further suppose that $f\in L^\infty_{I}$ is a positive function and 
\[ \theta := \sup_{(t,x)\in I\times \bbr^d} \sum_{l=1}^\infty \left(\int_I\int_{\bbr^d} G^{(l)}(t,x;s,y)\,\la(\dd s,\dd y)\right)^{1/p} < \infty. \]
Moreover, we assume that $v\in L^\infty_{I,\loc}$ is positive and satisfies
\beq\label{vineq} v(t,x)\leq f(t,x)+\left(\sum_{l=1}^\infty \int_I\int_{\bbr^d} G^{(l)}(t,x;s,y)(v(s,y))^{p\ga}\,\la(\dd s,\dd y)\right)^{1/p},\quad (t,x)\in I\times \bbr^d,\eeq
for some $\ga\in(0,1)$. Then $v\in L^\infty_{I}$ with $\|v\|_{L^\infty_{I}} \leq a$, where $a$ is the unique strictly positive solution to the equation $a-\|f\|_{L^\infty_{I}}- \theta a^\ga=0$.
\elem
\bpr The proof is a straightforward generalization of the arguments given in Example~\ref{inftyex}. We include it for the sake of completeness. Fix $T\in I$ and recall the definition of $\|\cdot\|_{G^{(l)},p}$ and $I_f(\cdot)$ from the proof of Lemma~\ref{comparison}. By \eqref{vineq}, it follows that
\[ \|v\|_{L^\infty_{I_T}} \leq \|I_f(v^\ga)\|_{L^\infty_{I_T}} \leq I_{\|f\|_{L^\infty_I}}(\|v\|_{L^\infty_{I_T}}^\ga). \]
By iteration of the last inequality, we deduce that $\|v\|_{L^\infty_{I_T}}\leq a_n(T)$ for all $n\in\bbn$ where $a_1(T):=\|v\|_{L^\infty_{I_T}}$ and $a_{n+1}(T)=I_{\|f\|_{L^\infty_I}}((a_n(T))^\ga) = \|f\|_{L^\infty_I} + \theta (a_n(T))^\ga$ for $n\in\bbn$. Straightforward analysis reveals that $\limsup_{n\to\infty} a_n(T)\leq a$, a number independent of $T$. Hence, $\|v\|_{L^\infty_I}\leq a$.
\epr

\section{Proof of the main theorems}\label{proofmain} 

\noindent \bff{Proof of Theorem~\ref{main1}.}\quad We show that Theorem~\ref{main1} is a special case of Theorem~\ref{main3}, or more precisely, that Assumption~\ref{Aass} is contained in Assumption~\ref{Cass}: setting $I=\bbr_+$ and $w\equiv1$ in Assumption~\ref{Cass}, it is not hard to see that the first six conditions break down to conditions (1)--(6) of Assumption~\ref{Aass}, and that condition (7) of Assumption~\ref{Cass} becomes superfluous. The only thing to show is that \eqref{partG0} implies \eqref{partG}. To this end, fix $T\in\bbr_+$, define
\[ \eps:=2^{-(p\vee1)}\left[(C_{\si,1}C_p^\mathrm{BDG})^p + C^p_{\si,1}\left(\sup_{(t,x)\in [0,T]\times\bbr^d} \int_0^t\int_{\bbr^d} |G(t,x;s,y)b_1(s,y)|\,\la(\dd s,\dd y)\right)^{p-1}\right]^{-1}, \]
and let $\mathcal{T}$ be a subdivision of $[0,T]$ such that \eqref{partG0} holds. Then we have for all $(t,x)\in [0,T]\times\bbr^d$ and $i=0,\ldots,k$ that 
\begin{align*} &~\sum_{l=1}^2 \left(\int_{t_i}^{t_{i+1}}\int_{\bbr^d}  G^{C,l}(t,x;s,y)\,\la(\dd s,\dd y)\right)^{1/(p\vee1)} \\
\leq&~ 2\left(\int_{t_i}^{t_{i+1}}\int_{\bbr^d}  G^{C,1}(t,x;s,y) + G^{C,2}(t,x;s,y)\,\la(\dd s,\dd y)\right)^{1/(p\vee1)} \\
\leq&~ \eps^{-1/(p\vee1)} \left(\int_{t_i}^{t_{i+1}}\int_{\bbr^d} G^A(t,x;s,y)\,\la(\dd s,\dd y)\right)^{1/(p\vee1)}<1,\end{align*}
which is \eqref{partG}. \halmos

\vspace{\baselineskip}
\noindent \bff{Proof of Corollary~\ref{cor1}.}\quad We check the conditions of Assumption~\ref{Aass}. (1), (2) and (3) are also assumed in the corollary. Regarding (4), (5) and (6), it is easy to see that because of \eqref{quasistat}, conditions \eqref{Gaussianfinite}, \eqref{jumpsfinite} and \eqref{meanstable0} split into separated conditions for both $G$ and $\La$, which are fulfilled thanks to \eqref{reg} and \eqref{localint}, respectively. Only (7) if left to be verified. Let $T\in\bbr_+$ be arbitrary and define $t_n^i:=i/n^2$ for $n\in\bbn$ and $i=0,\ldots,Tn^2$. Then, using the notation 
\beq\label{tilg} g^A:= \left(\int_\bbr |z|^p\,\pi_0(\dd z) + \|c\|_{L^\infty_{[0,T]}} \right) g^p + \|b_1\|_{L^\infty_{[0,T]}} g \bone_{\{p\geq1\}}, \eeq
we have for all $(t,x)\in[0,T]\times\bbr^d$ 
\[ \int_{t_n^i}^{t_n^{i+1}} \int_{\bbr^d} G^A(t,x;s,y) \,\dd (s,y) \leq \int_{t_n^i}^{t_n^{i+1}} \int_{\bbr^d} g^A(t-s,x-y) \,\dd(s,y) \leq \int_{(t-t_n^{i+1})\vee0}^{(t-t_n^i)\vee0} \int_{\bbr^d} g^A(s,y)\,\dd(s,y). \]
The right-hand side becomes arbitrarily small as $n\to\infty$, uniformly in $(t,x)\in[0,T]\times\bbr^d$ and $i=0,\ldots, Tn^2-1$. If not, there would exist some $\eps>0$ as well as for each $n\in\bbn$ some $\tau_n\in[0,T]$ and $i(n)\in\{0,\ldots,Tn^2-1\}$ such that
\[ \int_{(\tau_n-t_n^{i(n)+1})\vee0}^{(\tau_n-t_n^{i(n)})\vee0} \int_{\bbr^d} g^A(s,y)\,\dd(s,y) \geq \eps. \]
This, however, would contradict the dominated convergence theorem and the Borel-Cantelli lemma since $|((\tau_n-t_n^{i(n)})\vee0)-((\tau_n-t_n^{i(n)+1})\vee0)|\leq |t_n^{i(n)+1}-t_n^{i(n)}| = 1/n^2$. Thus, Corollary~\ref{cor1} is proved. \halmos

\vspace{\baselineskip}
\noindent \bff{Proof of Theorem~\ref{main2}.}\quad a) We first prove the existence of a solution to \eqref{SPDE-var}. To this end, define
\[ T_n := \inf\{t>0\colon |\La(\{t\}\times\bbr^d)|>n\},\quad n\in\bbn. \] 
Assumption~\ref{Bass}(3) implies that $(T_n)_{n\in\bbn}$ is a sequence of stopping times such that we have $T_n>0$ a.s. for each $n\in\bbn$ and $T_n\uparrow+\infty$ a.s. as $n\to\infty$. Next, we introduce for each $n\in\bbn$ a truncation of $\La$ in the following sense:
\begin{align*} \La^n(\dd t,\dd x) &:= B(\dd t,\dd x) + \La^\cc(\dd t,\dd x) + \int_\bbr z\bone_{\{|z|\leq1\}}\,(\mu-\nu)(\dd t,\dd x,\dd z) \\
&\quad + \int_\bbr z\bone_{\{1<|z|\leq n\}} \,\mu(\dd t,\dd x,\dd z).  \end{align*}
By Assumption~\ref{Bass}(4) we may assume without loss of generality that $Y_0\in B^q_{[0,\infty),\loc}$. Consequently, thanks to Assumption~\ref{Bass}(1) and (2) and Theorem~\ref{main1}, Equation~\eqref{SPDE-var} with $\La^n$ as driving noise has a unique solution $Y^n\in B^q_{[0,\infty),\loc}$. We claim that $Y:=Y^1 \bone_{\llbracket 0,T_1\rrbracket} + \sum_{n=2}^\infty Y^n \bone_{\rrbracket T_{n-1},T_n\rrbracket}$ is a solution to the original equation \eqref{SPDE-var} with $\La$. The predictability of $Y$ is clear. Now fix a (non-random) time $T\in\bbr_+$ and define
\[ \Om^n_T := \left\{\om\in\Om\colon \sup_{(t,x)\in [0,T]\times\bbr^d} |\La(\{(t,x)\})(\om)| \in [0,n]\right\}, \quad n\in\bbn. \]
By Assumption~\ref{Bass}(3) the sequence $(\Om^n_T)_{n\in\bbn}$ increases to $\Om$ up to a $\bbp$-null set. Moreover, we have $\bone_{\llbracket 0,T_k\rrbracket}(t) Y^k(t,x)=\bone_{\llbracket 0,T_k\rrbracket}(t) Y^n(t,x)$ a.s. for all $n\in\bbn$ and $k=1,\ldots,n$ as a consequence of the uniqueness statement of Theorem~\ref{main1} and the fact that $\bbp[T_k=t]=0$. Now part (1) of Theorem~\ref{main2} follows from the observation that for all $(t,x)\in[0,T]\times\bbr^d$ and $n\in\bbn$ we have a.s.
\begin{align*}
	&~\bone_{\Om^n_T}\int_0^t\int_{\bbr^d} G(t,x;s,y)\si(Y(s,y))\,\La(\dd s,\dd y) \\
	=&~ \bone_{\Om^n_T} \int_0^t\int_{\bbr^d} G(t,x;s,y)\left(\si(Y^1(s,y))\bone_{\llbracket 0, T_1 \rrbracket(s)} + \sum_{k=2}^n \si(Y^k(s,y))\bone_{\rrbracket T_{k-1},T_k\rrbracket}(s)\right) \,\La^n(\dd s,\dd y) \\
	=&~ \bone_{\Om^n_T} \int_0^t\int_{\bbr^d} G(t,x;s,y) \si(Y^n(s,y)) \,\La^n(\dd s,\dd y) = \bone_{\Om^n_T}Y^n(t,x) = \bone_{\Om^n_T} Y(t,x).
\end{align*}
To be utterly precise, for the transition from the second to the third line to be true, we must show that $J(\phi)$ and $J(\phi^\prime)$ as defined in \eqref{eqQ} are modifications of each other as soon as $\phi$ and $\phi^\prime$ are. But this follows from \eqref{momin2}. Finally, the uniqueness statement follows from that of Theorem~\ref{main1} by localization.

b) Next, we verify that the solution $Y$ found in a) belongs to $B^p_{[0,\infty),\loc}$ if also (5)--(10) of Assumption~\ref{Bass} are valid. We only carry out the proof for $p\geq1$. The case $p<1$ can be proved in the same fashion. Let $T\in\bbr_+$ and observe from a) that $Y$ equals $Y^n$ on $\Om^n_T$. Define $v^n(t,x):=\|Y^n(t,x)\|_{L^p}$ for $(t,x)\in[0,T]\times\bbr^d$, which is always finite because $Y^n\in B^q_{[0,\infty),\loc}$. Moreover, if we define $G^{D,l}$ as in \eqref{GD} with $w\equiv1$, then we have for all $(t,x)\in[0,T]\times\bbr^d$ according to Lemma~\ref{momin}(3) with $\rho=1$
\begin{align} \|Y(t,x)\bone_{\Om^n_T}\|_{L^p} \leq v^n(t,x) &\leq f(t,x) + \sum_{l=1}^4 C_{\si,2} \left(\int_0^t\int_{\bbr^d} G^{D,l}(t,x;s,y)(v^n(s,y))^p\,\la(\dd s,\dd y)\right)^{1/p}, \label{voltappear} \end{align}
where $f$ is the sum of the first three summands on the right-hand side of \eqref{momin3}. A priori, $G^{D,l}$ may depend on $n$ since it involves the underlying underlying L\'evy measure $\nu^n$. However, it is obvious that inequality \eqref{momin3} remains true if we use the original L\'evy measure $\nu$ to form $G^{D,l}$: the right-hand side of \eqref{voltappear} will only be enlarged. In this case, \eqref{voltappear} falls into the category of Lemma~\ref{comparison}(3). Indeed, Assumption~\ref{Bass}(10) guarantees that $f\in L^\infty_{[0,T]}$, and that the key assumption \eqref{partition} is met (note that the different constants appearing in $G^{D,l}$ compared to $G^B$ are irrelevant because $G^B$ satisfies the partition property \eqref{partG0} for all $\eps>0$). Thus, we have $v^n(t,x)\leq v(t,x)$ where $v\in L^\infty_{[0,T]}$ is again independent of $n$ and is the solution of the corresponding Volterra equation if we replace the second inequality sign in \eqref{voltappear} by equality. Taking the limit $n\to\infty$, we conclude 
\[ \|Y(t,x)\|_{L^p} = \lim_{n\to\infty} \|Y(t,x)\bone_{\Om^n_T}\|_{L^p} \leq v(t,x), \]
that is, $Y\in B^p_{[0,\infty),\loc}$. \halmos

\vspace{\baselineskip}
\noindent\bff{Proof of Corollary~\ref{cor2}.} \quad a) We begin with the first statement, for which we need to verify (2) and (3) of Assumption~\ref{Bass}. That (2) holds, follows from the proof of Corollary~\ref{cor1}, where we have shown that \eqref{quasistat}, \eqref{reg} and \eqref{localint} imply the validity of Assumption~\ref{Aass}(4)--(7). Notice that in the quasi-stationary case, it suffices to check Assumption~\ref{Bass}(2) only for $n=1$ because $\int_{1<|z|\leq n} |z|^q\,\pi_0(\dd z)$ is always finite and condition (3) of Corollary~\ref{cor2} is in force. That (3) of Assumption~\ref{Bass} holds, is due to \eqref{spacedecay}:
\[ \nu\big([0,T]\times\bbr^d\times[-1,1]^\Comp\big) \leq \int_0^T\int_{\bbr^d} \pi_1(t,x)\,\dd(t,x) \pi_0(|z|>1) < \infty. \]

b) For the second part we must prove (5)--(10) of Assumption~\ref{Bass}. (5) and (6) hold by hypothesis. Furthermore, since $p<q$  implies $|ab|^p_q\leq |a|^p_q |b|^q_p$ for all $a,b\in\bbr$, we have by \eqref{quasistat}
\begin{align*} &~\sup_{(t,x)\in[0,T]\times\bbr^d} \int_0^t\int_{\bbr^d}\int_{\bbr} |G(t,x;s,y)z|^p_q \,\nu(\dd s,\dd y,\dd z)\\
 \leq&~ \|\pi_1\|_{L^\infty_{[0,T]}} \int_\bbr |z|^p_q\,\pi_0(\dd z) \int_0^T\int_{\bbr^d} |g(t,x)|^q_p \,\dd(t,x) < \infty, \end{align*}
which implies (7) of Assumption~\ref{Bass}. Next, (8) is a direct consequence of condition (3) of the corollary. For (9) we choose $\al=q$ and $\beta=p$, which clearly satisfy (9c). For (9a) and (9b) first observe that 
\[ \left|b(t,x)+\int_\bbr z\bone_{\{|z|\in(1,A]\}}\,\pi(t,x,\dd z)\right| \leq \|b\|_{L^\infty_{[0,T]}} + \|\pi_1\|_{L^\infty_{[0,T]}} \int_{|z|>1} |z|^p\,\pi_0(\dd z) A^{1-p} \leq F_1 A^{1-p} \]
holds for all $A\in[1,\infty)$ if $F_1\in\bbr_+$ is chosen large enough. Second, if $q<1$, we have $b_0\equiv0$ by \eqref{reg}, which means that
\begin{align*} \left|b(t,x)-\int_\bbr z\bone_{\{|z|\in(a,1]\}}\,\pi(t,x,\dd z)\right| &= \left|\int_{\bbr} z\bone_{\{|z|\in(0,a]\}}\,\pi(t,x,\dd z)\right| \\
&\leq \|\pi_1\|_{L^\infty_{[0,T]}} \int_{|z|\leq1} |z|^q\,\pi_0(\dd z) a^{1-q}. \end{align*}
Finally, if $q\geq1$ we have
\[ \left|b(t,x)-\int_\bbr z\bone_{\{|z|\in(a,1]\}}\,\pi(t,x,\dd z)\right| \leq \|b\|_{L^\infty_{[0,T]}} + \|\pi_1\|_{L^\infty_{[0,T]}} \int_{|z|\leq1} |z|^q\,\pi_0(\dd z) a^{1-q} \leq F_0 a^{1-q} \] 
for all $a\in(0,1]$ and some constant $F_0\in\bbr_+$. Finally, condition (10) holds by the same arguments used in the proof of Corollary~\ref{cor1}. \halmos

\vspace{\baselineskip}
\noindent \bff{Proof of Theorem~\ref{main3}.}\quad We base the proof on a Picard iteration scheme, which parallels the construction of a solution to \eqref{volterra} in Lemma~\ref{comparison}. We define processes $Y^n\in\tilde\calp$ inductively as follows: starting with $Y^0(t,x) := Y_0(t,x)$, we assume that $Y^{n-1}\in B^{p,w}_{I,\loc}$ has already been constructed for some $n\in\bbn$. Define for each $(t,x)\in I\times\bbr^d$
\beq\label{iteration} Y^n(t,x) := Y_0(t,x) + \int_I \int_{\bbr^d} G(t,x;s,y)\si(Y^{n-1}(s,y))\,\La(\dd s, \dd y),\eeq
hereby choosing a predictable version of $Y^n$, cf. Lemma~\ref{predver}. Let $T\in I$. Then we have by Lemma~\ref{momin}(1) for all $(t,x)\in I_T\times\bbr^d$  
\begin{align*} &~\frac{\|Y^n(t,x)\|_{L^p}}{(w(t,x))^{1/(p\vee1)}} \leq \frac{\|Y_0(t,x)\|_{L^p}}{(w(t,x))^{1/(p\vee1)}} \\
&~\quad+ \sum_{l=1}^2 \left(\int_I \int_{\bbr^d}  \frac{G^{C,l}(t,x;s,y)}{C_{\si,1}^p}\left(\frac{|\si(0)|^{p\wedge1}+C_{\si,1}^{p\wedge1}\|Y^{n-1}(s,y)\|_{L^p}}{(w(s,y))^{1/(p\vee1)}}\right)^{p\vee1}\,\la(\dd s,\dd y)\right)^{1/(p\vee1)},
\end{align*}
which is finite by Assumption~\ref{Cass}. Thus, $Y^n\in B^{p,w}_{I,\loc}$ for all $n\in\bbn$. Next, Lemma~\ref{momin}(2) implies that $u^n:=Y^n - Y^{n-1}$ satisfies 
	\beq \label{unit} \frac{\|u^{n+1}(t,x)\|_{L^p}}{(w(t,x))^{1/(p\vee1)}}\leq \sum_{l=1}^2 \left(\int_I \int_{\bbr^d}  G^{C,l}(t,x;s,y)\left(\frac{\|u^n(s,y)\|_{L^p}}{(w(s,y))^{1/(p\vee1)}}\right)^{p\vee1}\,\la(\dd s,\dd y)\right)^{1/(p\vee1)}\eeq
for all $(t,x)\in I\times\bbr^d$, which is a recursive relation as in Lemma~\ref{comparison}(1). Note that the key hypothesis \eqref{partition} is fulfilled because of Assumption~\ref{Cass}(8). We conclude that $\sum_{n=1}^\infty \|u^n\|_{B^{p,w}_{I_T}}<\infty$, in other words, $Y^n$ converges in $B^{p,w}_{I_T}$ to some limit $Y$. Applying Lemma~\ref{momin}(2) to $\phi_1:=Y$ and $\phi_2:=Y^{n-1}$, the convergence $Y^{n-1} \to Y$ also implies that $J(Y^{n-1})=Y^n \to J(Y)$ in $B^{p,w}_{I_T}$, that is, $Y$ indeed satisfies \eqref{SPDE-var}. The uniqueness of the solution to \eqref{SPDE-var} follows if we substitute $u^n$ in \eqref{unit} by the difference of two solutions. Since $T\in I$ is arbitrary, Theorem~\ref{main3} follows. \halmos

\vspace{\baselineskip}
\noindent \bff{Proof of Corollary~\ref{cor3}.}\quad We verify Assumption~\ref{Cass} for $I=\bbr$ and $w\equiv1$. (1), (2) and (3) hold by hypothesis; (4), (5) and (6) are consequences of \eqref{quasistat2}, \eqref{reg2}, \eqref{intconst1} and \eqref{intconst2}. Moreover, condition (7) of Assumption~\ref{Cass} is redundant such that it remains to verify (8). To this end, define
\begin{align*} g^{C,1}&:= (C_{\si,1}C^\mathrm{BDG}_p)^p(\zeta_p + \|c\|_{L^\infty_\bbr})g^p,\\
g^{C,2}&:=(C_{\si,1}\|b_1\|_{L^\infty_\bbr})^p\left(\int_0^\infty\int_{\bbr^d} g(t,x)\,\dd(t,x)\right)^{p-1}g\bone_{\{p\geq1\}}.
\end{align*}
Then, for any subdivision $\mathcal{T}\colon -\infty=t_0 < \ldots < t_{k+1} = T$, all $(t,x)\in(-\infty,T]\times\bbr^d$ and $i=0,\ldots,k$, we have by \eqref{intconst1} and \eqref{intconst2}
\begin{align*} \sum_{l=1}^2 \left(\int_{t_i}^{t_{i+1}}\int_{\bbr^d} G^{C,l}(t,x;s,y) \,\dd (s,y)\right)^{1/(p\vee1)} &\leq \sum_{l=1}^2\left(\int_{-\infty}^t\int_{\bbr^d} g^{C,l}(t-s,x-y)\,\dd(s,y)\right)^{1/(p\vee1)} \\
&= \sum_{l=1}^2 \left( \int_0^\infty \int_{\bbr^d} g^{C,l}(t,x)\,\dd(t,x) \right)^{1/(p\vee1)} < 1. \end{align*}
\halmos

\vspace{\baselineskip}
\noindent \bff{Proof of Theorem~\ref{further}.}\quad a) Fix $T\in I$ and choose $(t,x), (\tau,\xi)\in I_T\times\bbr^d$. Then similar calculations as in Lemma~\ref{momin}(2) lead to
\begin{align*} \|Y(t,x)-Y(\tau,\xi)\|_{L^p} &\leq \sum_{l=1}^2 \left(\int_I\int_{\bbr^d} \tilde G^{(l)}(t,x;\tau,\xi;s,y) \left(\frac{\|\si(Y(s,y))\|_{L^p}}{(w(s,y))^{1/(p\vee1)}}\right)^{p\vee1}\,\la(\dd s,\dd y)\right)^{1/(p\vee1)}\\
	&\leq \|\si(Y)\|_{B^{p,w}_{I_T}} \sum_{l=1}^2 \left(\int_I\int_{\bbr^d} \tilde G^{(l)}(t,x;\tau,\xi;s,y) \,\la(\dd s,\dd y)\right)^{1/(p\vee1)},
\end{align*}
where 
\begin{align*}
\tilde G^{(1)}(t,x;\tau,\xi;s,y) &:= (C^\mathrm{BDG}_p)^p |G(t,x;s,y)-G(\tau,\xi;s,y)|^p \left(\int_\bbr |z|^p \,\pi(s,y,\dd z) + c(s,y)\right)w(s,y),\\
\tilde G^{(2)}(t,x;\tau,\xi;s,y) &:= \left(\int_I\int_{\bbr^d} |[G(t,x;s,y)-G(\tau,\xi;s,y)]b_1(s,y)|\,\la(\dd s,\dd y)\right)^{p-1}\\
&\quad\cdot|[G(t,x;s,y)-G(\tau,\xi;s,y)]b_1(s,y)|w(s,y)\bone_{\{p\geq1\}}.
\end{align*} 
The claim now follows from \eqref{Lpcontcond} because Assumption~\ref{Cass}(7) implies
\begin{align*} &~\sup_{(t,x),(\tau,\xi)\in I_T\times\bbr^d} \left(\int_I\int_{\bbr^d} |[G(t,x;s,y)-G(\tau,\xi;s,y)]b_1(s,y)|\,\la(\dd s,\dd y)\right)^{p-1} \\
\leq&~2\sup_{(t,x)\in I_T\times\bbr^d} \left(\int_I\int_{\bbr^d} |G(t,x;s,y)b_1(s,y)|\,\la(\dd s,\dd y)\right)^{p-1}<\infty. \end{align*}

b) In the situation of Corollary~\ref{cor3} with $G$ in convolution form, we have 
\begin{align*} \int_\bbr\int_{\bbr^d} \tilde G(t,x;\tau,\xi;s,y) \,\dd (s,y)&\leq(\zeta_p + \|c\|_{L^\infty_\bbr}) \int_{\bbr}\int_{\bbr^d} \big|g(t-s,x-y)-g(\tau-s,\xi-y)\big|^p\,\dd(s,y)\\
&\quad+ \|b_1\|_{L^\infty_\bbr} \bone_{\{p\geq1\}} \int_\bbr\int_{\bbr^d} |g(t-s,x-y)-g(\tau-s,\xi-y)|\,\dd(s,y) \\
 &= (\zeta_p + \|c\|_{L^\infty_\bbr}) \int_{\bbr}\int_{\bbr^d} \big|g(s+h,y+\eta)-g(s,y)\big|^p\,\dd(s,y) \\
 &\quad+ \|b_1\|_{L^\infty_\bbr} \bone_{\{p\geq1\}} \int_\bbr\int_{\bbr^d} |g(s+h,y+\eta)-g(s,y)\big|\,\dd(s,y) \to 0 \end{align*}
because $(h,\eta)=(|t-\tau|,|x-\xi|)\to 0$, cf. \citep[Lemma~0.12]{Folland95}.

c) Let $T\in I$, $\bar p:=p\vee1$ and define $v(t,x):=w^{-1/\bar p}(t,x)\|Y(t,x)-Y^\prime(t,x)\|_{L^p}$ as well as $v_0(t,x):=w^{-1/\bar p}(t,x)\|Y_0(t,x)-Y^\prime_0(t,x)\|_{L^p}$. Furthermore, choose $k\in\bbn$ and a partition $I_T=I_1\cup\ldots\cup I_k$ such that \eqref{partG} is satisfied. Next, recall from \eqref{Gnorms} the definition of $\|\phi\|_{G^{(l)},\bar p}(t,x)$ and $\|\phi\|_{G^{(l)},\bar p,j}(t,x)$ for $(t,x)\in I_T\times\bbr^d$, $l=1,2$ and $j=1,\ldots,k$. From Lemma~\ref{momin}(2) we deduce 
\beq\label{v-it} v \leq v_0 + \sum_{l=1}^2 \|v\|_{G^{(l)},\bar p} \leq v_0 + \sum_{j=1}^k \sum_{l=1}^2 \|v\|_{G^{(l)},\bar p,j}. \eeq
By the same arguments as in the proof of Lemma~\ref{comparison}(1), iterating \eqref{v-it} $N$ times produces
\[ \|v\|_{L^\infty_{I_T}}\leq \|v_0\|_{L^\infty_{I_T}} \sum_{n=0}^{N-1} \binom{n+k-1}{n}\rho^n + \|v\|_{L^\infty_{I_T}}\binom{N+k-1}{N}\rho^N,  \]
with $\rho<1$ being the left-hand side of \eqref{partG}. Letting $N\to\infty$ leads to the assertion.\halmos

\vspace{\baselineskip}
\noindent\bff{Proof of Theorem~\ref{statsol}.}\quad It suffices to prove the case where \eqref{cag0} holds. Since $Y\in B^p_{\bbr,\loc}$ is constructed as the limit of the Picard iterates $Y^n$ in \eqref{iteration}, it suffices to prove that $Y^n$, $Y_0$ and $\La$ are jointly stationary for all $n\in\bbn$. By induction, we assume that $Y^{n-1}$ is jointly stationary with $\La$ and $Y_0$ (that $Y_0$ is, holds by assumption). First, we assume that $g$ is bounded and has compact support in $\bbr_+\times\bbr^d$, which obviously implies that \eqref{supbounded} holds for arbitrary $\eps>0$. Moreover, $Y^{n-1}$ is $L^p$-continuous because $Y^0$ is by hypothesis and thus also $Y^{n-1}$ for general $n$ by the same arguments as in the proof of Theorem~\ref{further}(2). Next, we fix $(t,x),(h,\eta)\in \bbr\times\bbr^d$ and define for $N\in\bbn$ and $i=0,\ldots,N^2$ the time points $s_i^N:=t-N+i/N$. Moreover, we set $Q_N:=\big(0,(1/N,\ldots,1/N)\big]$ and $\Gamma_N:=\big\{(i_1/N,\ldots,i_d/N)\colon i_1,\ldots,i_d\in\{-N^2,\ldots,N^2\}\big\}$. Lemma~\ref{Riemann} now gives
\allowdisplaybreaks
\begin{align*} &~Y^n(t+h,x+\eta)=Y_0(t+h,x+\eta) + \int_{-\infty}^{t+h}\int_{\bbr^d} g(t+h-s,x+\eta-y)\si(Y^{n-1}(s,y)\,\La(\dd s,\dd y)\\
=&~Y_0(t+h,x+\eta) + \int_{-\infty}^t \int_{\bbr^d} g(t-s,x-y)\si(Y^{n-1}(s+h,y+\eta)\,\La(h+\dd s,\eta+\dd y)\\
=&~Y_0(t+h,x+\eta) + L^p\!\!-\!\!\lim_{N\to\infty} \sum_{i=0}^{N^2-1}\sum_{y_j^N \in \Gamma_N} g(t-s_i^N,x-y_j^N)\si(Y^{n-1}(s_i^N+h,y_j^N+\eta))\\
&~\quad\cdot\La\big((s_i^N+h,s_{i+1}^N+h) \times (y_j^N + \eta + Q_N) \big)\\
\eqd&~Y_0(t,x) + L^p\!\!-\!\!\lim_{N\to\infty} \sum_{i=0}^{N^2-1}\sum_{y_j^N \in \Gamma_N} g(t-s_i^N,x-y_j^N)\si(Y^{n-1}(s_i^N,y_j^N))\La\big((s_i^N,s_{i+1}^N) \times (y_j^N + Q_N) \big)\\
=&~Y^n(t,x). \end{align*} 
\allowdisplaybreaks[0]
The calculation remains valid when we consider joint distributions with $Y_0$ and $\La$, and when we extend it to $n$ space--time points. So the theorem is proved for bounded functions $g$ with compact support. For general functions $g$ we notice that property \eqref{cag0} implies that we can write $g=\sum_{i=1}^\infty g_i$ where each $g_i$ is bounded with compact support. The theorem follows since the calculation above is invariant under summation and taking limits. \halmos

\vspace{\baselineskip}
\noindent\bff{Proof of Theorem~\ref{assstat}.}\quad Let $Y\in B^{p,w}_{I,\loc}$ be a solution to \eqref{SPDE-var}. Then we have $v\in L^\infty_{I,\loc}$ where $v$ is defined by $v(t,x):=w^{-{1/(p\vee1)}}(t,x)\|Y(t,x)\|_{L^p}$. The claim is that $v$ also belongs to $L^\infty_I$. We only consider the case $p\in[1,2]$, the case $p\in(0,1)$ can be treated analogously. First, we suppose that Assumption~\ref{Dass}(6a) holds. In this case, it follows from Lemma~\ref{momin}(3) that there exists some $\rho\in(0,1)$ with
\beq\label{eq1} v(t,x) \leq f(t,x) + \sum_{l=1}^4 C_{\si,2} \left(\int_I \int_{\bbr^d} G^{D,l}(t,x;s,y)(w(s,y))^{\rho-1} (v(s,y))^{p\rho} \,\la(\dd s,\dd y)\right)^{1/p}, \eeq
where $f$ denotes the sum of the first three terms on the right-hand side of \eqref{momin3}. By hypothesis, the functions $w^{-1}$, $w^{-1/p}$ and $w^{\rho-1}$ are uniformly bounded on $I\times\bbr^d$, which means that $f$ belongs to $L^\infty_I$. Consequently, Lemma~\ref{assympholder} together with (3), (4) and (5) of Assumption~\ref{Dass} shows that $v\in L^\infty_I$. Now suppose that Assumption~\ref{Dass}(6b) holds. Then, by replacing $r$ in \eqref{eq1} by $1$, the claim follows from Lemma~\ref{comparison}(3) and assumption \eqref{D-6b}.\halmos

\subsection*{Acknowledgement}
I take pleasure in thanking Claudia Kl\"uppelberg and Jean Jacod for their valuable advice and careful proofreading. Support from the graduate program TopMath at Technische Universit\"at M\"unchen and the Studienstiftung des deutschen Volkes is gratefully acknowledged.

\addcontentsline{toc}{section}{References}
\bibliographystyle{plainnat}
\bibliography{bib-SPDE}

\begin{thebibliography}{32}
\providecommand{\natexlab}[1]{#1}
\providecommand{\url}[1]{\texttt{#1}}
\expandafter\ifx\csname urlstyle\endcsname\relax
  \providecommand{\doi}[1]{doi: #1}\else
  \providecommand{\doi}{doi: \begingroup \urlstyle{rm}\Url}\fi

\bibitem[A{\"i}t-Sahalia and Jacod(2014)]{AitSahalia14}
Y.~A{\"i}t-Sahalia and J.~Jacod.
\newblock \emph{High-Frequency Financial Econometrics}.
\newblock Princeton University Press, Princeton, 2014.

\bibitem[Albeverio et~al.(1998)Albeverio, Wu, and Zhang]{Albeverio98}
S.~Albeverio, J.-L. Wu, and T.-S. Zhang.
\newblock Parabolic {SPDEs} driven by {P}oisson white noise.
\newblock \emph{Stoch. Process. Appl.}, 74\penalty0 (1):\penalty0 21--36, 1998.

\bibitem[Applebaum and Wu(2000)]{Applebaum00}
D.~Applebaum and J.-L. Wu.
\newblock Stochastic partial differential equations driven by {L}{\'e}vy
  space--time white noise.
\newblock \emph{Random Oper. Stoch. Equ.}, 8\penalty0 (3):\penalty0 245--259,
  2000.

\bibitem[Asmussen(2003)]{Asmussen03}
S.~Asmussen.
\newblock \emph{Applied Probability and Queues}.
\newblock Springer, New York, 2nd edition, 2003.

\bibitem[Balan(2014)]{Balan14}
R.M. Balan.
\newblock {SPDE}s with {$\alpha$}-stable {L}{\'e}vy noise: a random field
  approach.
\newblock \emph{Int. J. Stoch. Anal.}, 2014.
\newblock Article ID 793275, 22 pages.

\bibitem[Barndorff-Nielsen and Schmiegel(2004)]{BN04}
O.E. Barndorff-Nielsen and J.~Schmiegel.
\newblock L{\'e}vy-based spatial--temporal modelling, with applications to
  turbulence.
\newblock \emph{Russ. Math. Surv.}, 59\penalty0 (1):\penalty0 65--90, 2004.

\bibitem[Barndorff-Nielsen et~al.(2011)Barndorff-Nielsen, Benth, and
  Veraart]{BN11-2}
O.E. Barndorff-Nielsen, F.E. Benth, and A.E.D. Veraart.
\newblock Ambit processes and stochastic partial differential equations.
\newblock In G.~Di Nunno and B.~{\O}ksendal, editors, \emph{Advanced
  Mathematical Methods for Finance}, pages 35--74. Springer, Berlin, 2011.

\bibitem[Barndorff-Nielsen et~al.(2015)Barndorff-Nielsen, Benth, and
  Veraart]{BN15}
O.E. Barndorff-Nielsen, F.E. Benth, and A.E.D. Veraart.
\newblock Recent advances in ambit stochastics with a view towards
  tempo-spatial stochastic volatility/intermittency.
\newblock \emph{Banach Cent. Publ.}, 104:\penalty0 25--60, 2015.

\bibitem[Bichteler(2002)]{Bichteler02}
K.~Bichteler.
\newblock \emph{Stochastic Integration with Jumps}.
\newblock Cambridge University Press, Cambridge, 2002.

\bibitem[Bichteler and Jacod(1983)]{Bichteler83}
K.~Bichteler and J.~Jacod.
\newblock Random measures and stochastic integration.
\newblock In G.~Kallianpur, editor, \emph{Theory and Application of Random
  Fields}, pages 1--18. Springer, Berlin, 1983.

\bibitem[Chong and Kl\"uppelberg(2015)]{Chong14}
C.~Chong and C.~Kl\"uppelberg.
\newblock Integrability conditions for space--time stochastic integrals: Theory
  and applications.
\newblock \emph{Bernoulli}, 21\penalty0 (4):\penalty0 2190--2216, 2015.

\bibitem[Couchran et~al.(1995)Couchran, Lee, and Potthoff]{Couchran95}
W.G. Couchran, J.-S. Lee, and J.~Potthoff.
\newblock Stochastic {V}olterra equations with singular kernels.
\newblock \emph{Stoch. Process. Appl.}, 56\penalty0 (2):\penalty0 337--349,
  1995.

\bibitem[Coutin and Decreusefond(2001)]{Coutin01}
L.~Coutin and L.~Decreusefond.
\newblock Stochastic {V}olterra equations with singular kernels.
\newblock In A.B. Cruzeiro and J.-C. Zambrini, editors, \emph{Stochastic
  Analysis and Mathematical Physics}, pages 39--50. Birkh{\"a}user, Boston,
  2001.

\bibitem[Dalang(1999)]{Dalang99}
R.C. Dalang.
\newblock Extending martingale measures stochastic integral with applications
  to spatially homogeneous {S.P.D.E's}.
\newblock \emph{Electron. J. Probab.}, 4\penalty0 (6), 1999.
\newblock 24 pages.

\bibitem[Dalang and Quer-Sardanyons(2011)]{Dalang11}
R.C. Dalang and L.~Quer-Sardanyons.
\newblock Stochastic integrals for spde's: {A} comparison.
\newblock \emph{Expo. Math.}, 29\penalty0 (1):\penalty0 67--109, 2011.

\bibitem[Folland(1995)]{Folland95}
G.B. Folland.
\newblock \emph{Introduction to Partial Differential Equations}.
\newblock Princeton University Press, Princeton, 2nd edition, 1995.

\bibitem[Gripenberg et~al.(1990)Gripenberg, Londen, and Steffans]{Gripenberg90}
G.~Gripenberg, S.-O. Londen, and O.~Steffans.
\newblock \emph{Volterra Integral and Functional Equations}.
\newblock Cambridge University Press, Cambridge, 1990.

\bibitem[Jacod and Shiryaev(2003)]{Jacod03}
J.~Jacod and A.N. Shiryaev.
\newblock \emph{Limit Theorems for Stochastic Processes}.
\newblock Springer, Berlin, 2nd edition, 2003.

\bibitem[Kov{\'a}cs et~al.(2015)Kov{\'a}cs, Lindner, and Schilling]{Kovacs15}
M.~Kov{\'a}cs, F.~Lindner, and R.L. Schilling.
\newblock Weak convergence of finite element approximations of linear
  stochastic evolution equations with additive l{\'e}vy noise.
\newblock Preprint under arXiv:1411.1051 [math.PR], 2015.

\bibitem[Lebedev(1995)]{Lebedev95}
V.A. Lebedev.
\newblock The {F}ubini theorem for stochastic integrals with respect to ${L}^0
  $-valued random measures depending on a parameter.
\newblock \emph{Theory Probab. Appl.}, 40\penalty0 (2):\penalty0 285--293,
  1995.

\bibitem[Mueller(1998)]{Mueller98}
C.~Mueller.
\newblock The heat equation with {L}{\'e}vy noise.
\newblock \emph{Stoch. Process. Appl.}, 74\penalty0 (1):\penalty0 67--82, 1998.

\bibitem[Mytnik(2002)]{Mytnik02}
L.~Mytnik.
\newblock Stochastic partial differential equation driven by stable noise.
\newblock \emph{Probab. Theory Relat. Fields}, 123\penalty0 (2):\penalty0
  157--201, 2002.

\bibitem[Os{\k{e}}kowski(2012)]{Osekowski12}
A.~Os{\k{e}}kowski.
\newblock \emph{Sharp Martingale and Semimartingale Inequalities}.
\newblock Birkh{\"a}user, Basel, 2012.

\bibitem[Peszat and Zabczyk(2007)]{Peszat07}
S.~Peszat and J.~Zabczyk.
\newblock \emph{Stochastic Partial Differential Equations with L{\'e}vy Noise}.
\newblock Cambridge University Press, Cambridge, 2007.

\bibitem[Podolskij(2015)]{Podolskij14}
M.~Podolskij.
\newblock Ambit fields: Survey and new challenges.
\newblock In R.H. Mena, J.C. Pardo, V.~Rivero, and G.U. Bravo, editors,
  \emph{XI Symposium on Probability and Stochastic Processes}, pages 241--279.
  Springer, Cham, 2015.

\bibitem[Protter(1985)]{Protter85}
P.~Protter.
\newblock Volterra equations driven by semimartingales.
\newblock \emph{Ann. Probab.}, 13\penalty0 (2):\penalty0 519--530, 1985.

\bibitem[Rajput and Rosi{\'n}ski(1989)]{Rajput89}
B.S. Rajput and J.~Rosi{\'n}ski.
\newblock Spectral representations of infinitely divisible processes.
\newblock \emph{Probab. Theory Relat. Fields}, 82\penalty0 (3):\penalty0
  451--487, 1989.

\bibitem[Rei{\ss} et~al.(2006)Rei{\ss}, Riedle, and van Gaans]{Reiss06}
M.~Rei{\ss}, M.~Riedle, and O.~van Gaans.
\newblock Delay differential equations driven by {L}{\'e}vy processes:
  {S}tationarity and {F}eller properties.
\newblock \emph{Stoch. Process. Appl.}, 116\penalty0 (10):\penalty0 1409--1432,
  2006.

\bibitem[{Saint Loubert Bi{\'e}}(1998)]{SLB98}
E.~{Saint Loubert Bi{\'e}}.
\newblock {\'E}tude d'une {EDPS} conduite par un bruit poissonnien.
\newblock \emph{Probab. Theory Relat. Fields}, 111\penalty0 (2):\penalty0
  287--321, 1998.

\bibitem[Stricker and Yor(1978)]{Stricker78}
C.~Stricker and M.~Yor.
\newblock Calcul stochastique d{\'e}pendant d'un param{\`e}tre.
\newblock \emph{Z. Wahrscheinlichkeitstheorie verw. Geb.}, 45\penalty0
  (2):\penalty0 109--133, 1978.

\bibitem[Walsh(1986)]{Walsh86}
J.B. Walsh.
\newblock An introduction to stochastic partial differential equations.
\newblock In P.L. Hennequin, editor, \emph{\'Ecole d'\'Et\'e de Probabilit\'es
  de Saint Flour XIV - 1984}, volume 1180 of \emph{Lecture Notes in
  Mathematics}, pages 265--439. Springer, Berlin, 1986.

\bibitem[Wang(2008)]{Wang08}
Z.~Wang.
\newblock Existence and uniqueness of solutions to stochastic {V}olterra
  equations with singular kernels and non-{L}ipschitz coefficients.
\newblock \emph{Stat. Probab. Lett.}, 78\penalty0 (9):\penalty0 1062--1071,
  2008.

\end{thebibliography}
\end{document}